\documentclass[11pt,a4paper]{article}
\usepackage{a4wide}
\usepackage[utf8]{inputenc}
\usepackage{latexsym}
\usepackage{amsmath}
\usepackage{amsthm}
\usepackage{amssymb,enumerate}
\usepackage[usenames]{color}
\usepackage{graphicx}
\usepackage{fancybox}
\usepackage[normalem]{ulem}
\usepackage{hyperref}
\usepackage{pgfplots}
\pgfplotsset{compat=1.15}
\usepackage{mathrsfs}
\usetikzlibrary{arrows}
\definecolor{qqqqff}{rgb}{0.,0.,1.}
\definecolor{ududff}{rgb}{0.30196078431372547,0.30196078431372547,1.}
\definecolor{ffqqqq}{rgb}{1.,0.,0.}
\definecolor{xdxdff}{rgb}{0.49019607843137253,0.49019607843137253,1.}
\definecolor{verdeosc}{rgb}{0,0.5,0} 
\theoremstyle{plain}
\newtheorem{proposition}{Proposition}[section]
\newtheorem{theorem}[proposition]{Theorem}
\newtheorem{lemma}[proposition]{Lemma}
\newtheorem{corollary}[proposition]{Corollary}

\theoremstyle{definition}
\newtheorem{definition}[proposition]{Definition}

\newtheorem{remark}[proposition]{Remark}

\usepackage{geometry}
\usepackage{subeqnarray}
\usepackage{enumerate}
\usepackage{color}

\makeatletter
\DeclareRobustCommand\widecheck[1]{{\mathpalette\@widecheck{#1}}}
\def\@widecheck#1#2{%
	\setbox\z@\hbox{\m@th$#1#2$}%
	\setbox\tw@\hbox{\m@th$#1%
		\widehat{%
			\vrule\@width\z@\@height\ht\z@
			\vrule\@height\z@\@width\wd\z@}$}%
	\dp\tw@-\ht\z@
	\@tempdima\ht\z@ \advance\@tempdima2\ht\tw@ \divide\@tempdima\thr@@
	\setbox\tw@\hbox{%
		 \raise\@tempdima\hbox{\scalebox{1}[-1]{\lower\@tempdima\box
				\tw@}}}%
	{\ooalign{\box\tw@ \cr \box\z@}}}
\makeatother
\renewcommand{\check}[1]{\widecheck{#1}}

\usepackage{bbm}

\newcommand{\RR}{\mathbb{R}}
\newcommand{\CC}{\mathbb{C}}
\newcommand{\NN}{\mathbb{N}}

\newcommand{\M}{\mathbb{M}}
\newcommand{\m}{\boldsymbol{m}}
\newcommand{\LL}{\mathbb{L}}
\newcommand{\W}{\mathbb{W}}
\newcommand{\G}{\mathbb{G}}

\let\on=\operatorname

\newenvironment{demo}[1]{\par\smallskip\noindent{\bf #1.}}{\par\smallskip}

\makeatletter
\@namedef{subjclassname@2020}{%
	\textup{2020} Mathematics Subject Classification}
\makeatother
\def\a{\alpha}

\definecolor{azulosc}{rgb}{0.2,0.1,0.7}
\definecolor{granate}{rgb}{0.6,0,0.3} %col4
\definecolor{verdeosc}{RGB}{46, 139, 87} %colororiginal
\definecolor{rojo}{RGB}{219,0,0}                %col9

\makeatletter
\DeclareRobustCommand\widecheck[1]{{\mathpalette\@widecheck{#1}}}
\def\@widecheck#1#2{%
    \setbox\z@\hbox{\m@th$#1#2$}%
    \setbox\tw@\hbox{\m@th$#1%
       \widehat{%
          \vrule\@width\z@\@height\ht\z@
          \vrule\@height\z@\@width\wd\z@}$}%
    \dp\tw@-\ht\z@
    \@tempdima\ht\z@ \advance\@tempdima2\ht\tw@ \divide\@tempdima\thr@@
    \setbox\tw@\hbox{%
       \raise\@tempdima\hbox{\scalebox{1}[-1]{\lower\@tempdima\box
\tw@}}}%
    {\ooalign{\box\tw@ \cr \box\z@}}}
\makeatother
\begin{document}

\title{Stability properties of ultraholomorphic classes of Roumieu-type defined by weight matrices}
\author{Javier Jiménez-Garrido \and Ignacio Miguel-Cantero \and Javier Sanz \and Gerhard Schindl}
\date{\today}

\maketitle

\begin{abstract}
We characterize several stability properties, such as inverse or composition closedness, for ultraholomorphic function classes of Roumieu type defined in terms of a weight matrix. In this way we transfer and extend known results from J. Siddiqi and M. Ider, from the weight sequence setting and in sectors not wider than a half-plane, to the weight matrix framework and for sectors in the Riemann surface of the logarithm with arbitrary opening. The key argument rests on the construction, under suitable hypotheses, of characteristic functions in these classes for unrestricted sectors. As a by-product, we obtain new stability results when the growth control in these classes is expressed in terms of a weight sequence, or of a weight function in the sense of Braun-Meise-Taylor.
\par\medskip

 \end{abstract}

\section{Introduction}

When dealing with function spaces (usually called classes) it is very interesting to decide whether the usual operations (pointwise product, composition, algebraic inversion, differentiation, integration, etc.) on the functions of the space provide new functions inside it. These stability properties play a crucial role in the setting and the solution of, for example, algebraic, differential or integro-differential equations in the class.

In the literature one can frequently find the so-called ultradifferentiable classes, both in the Carleman and the Braun-Meise-Taylor sense, whose elements are smooth functions defined on open subsets of $\RR^n$ (or possibly germs at a point) such that the rate of growth of their successive derivatives is controlled (except for a geometric factor) in terms of a given sequence of positive real numbers in the first case, or of a given weight function in the second one. Moreover, depending on the choice of a universal or existential quantifier for the geometric factor in the estimates, one can consider Beurling- or Roumieu-like classes in both situations. The study of stability under inversion (or division) in these frameworks has a long history, see the works of Rudin~\cite{RudinDivision}, Bruna~\cite{BrunaInverse} and Siddiqi~\cite{Siddiqi}, and also composition has been studied in Fern\'andez and Galbis~\cite{FernandezGalbis06}. Recently, the introduction of classes associated with a weight matrix, by the fourth author of this paper~\cite{diploma,dissertation}, which strictly encompass those classes mentioned before, has led him and Rainer~\cite{compositionpaper,characterizationstabilitypaper}
to the characterization of stability under different operations in terms of conditions for the weight matrix under consideration, so giving a satisfactory general solution to these problems.

In connection with the asymptotic theory of solutions for differential and difference equations around singular points in the complex domain, it is natural to consider the complex analogue of such classes, usually called ultraholomorphic classes. They consist of holomorphic functions in sectorial regions in the Riemann surface of the logarithm (the singular point is assumed to be at 0, the vertex of the region) whose derivatives admit again suitable estimates of Roumieu type in terms of a sequence of positive real numbers, which in the applications is typically a Gevrey sequence $(p!^{a})_{p\in\NN_0}$ for some $a>1$. The study of stability properties in such classes is well-known for the Gevrey ones, see~\cite{BalserUTX}, but already in 1987 Ider and Siddiqi~\cite{IderSiddiqi} studied stability under composition with analytic functions and under inversion for general Carleman-Roumieu classes in unbounded sectors not wider than a half-plane. Our aim is to extend their results in several senses: (1) we consider Roumieu classes defined by weight matrices, so including in our considerations those of Carleman type and those defined by a weight function, as in the ultradifferentiable setting; (2) we are able to deal with classes defined in sectors of arbitrary opening in the Riemann surface of the logarithm, and (3) we extend the list of stability properties, including that of composition closedness.
It is important to note that, in the case of classes given by a weight function, a fundamental role in the stability properties is played by the condition that this function is equivalent to a concave weight function, what amounts to the root almost increasing property for the associated weight matrix.

The main novelties arise from two different sources. On the one hand, the techniques coming with the weight matrix structure allow for a better understanding of the conditions usually appearing in such stability results, and provide a clear way to establish results for the weight sequence and weight function approach. Indeed, our results extend the known ones for Carleman classes, and they match, in the limit when the opening of the sector tends to 0, with the ones for ultradifferentiable classes on a half-line. On the other hand, the main statements heavily rest on the construction of so-called characteristic functions in Carleman-Roumieu ultraholomorphic classes in sectors of arbitrary opening. These functions are those in a class which cannot belong to a class strictly contained in the original one, and so are in a sense maximal within the class. While Ider and Siddiqi only got such functions in narrow sectors, the work of Rodríguez Salinas~\cite{salinas62} provides indeed the key facts for working in general sectors, and this is in turn crucial for our purposes.

The paper is organized as follows. Section~\ref{sect.WeightCondit} contains all the preliminary information about sequences, weight functions and weight matrices. For the ultraholomorphic classes introduced in Section~\ref{sect.UltraholomClass} we show how to construct characteristic functions in Section~\ref{sect.CharactFunctions}. The stability results for classes associated with weight matrices are given in Section~\ref{sect.StabilityPropert}, and Section~\ref{sect.StabilityWeightFUnctCase} is devoted to their particularization to the case of classes induced by a weight function.
{Finally, we present in Section~\ref{sect.Examples} some examples, including those of Gevrey and $q$-Gevrey classes, in order to illustrate the obtained results.}

\section{Preliminaries on sequences, weight functions and weight matrices}\label{sect.WeightCondit}
%\subsection{General notation}
 %With $[\cdot]$ we denote either the Roumieu case $\{\cdot\}$ or the Beurling case $(\cdot)$ but not mixing the cases (for function spaces, classes of weighted sequences and conditions).

\subsection{Weight sequences}\label{ssectweightsequences}

We write $\NN_0:=\{0,1,2,\dots\}$ and $\NN:=\{1,2,3,\dots\}$. In what follows, we always denote by $\M=(M_j)_j\in\RR_{>0}^{\NN_0}$ a sequence with $M_0=1$,  we also use $\check{\M}=(\check{M}_j)_j$ defined by $\check{M}_j:=\frac{M_j}{j!}$ and the sequence of quotients associated $\m=(m_j)_j$ defined by $m_j:=M_{j+1}/M_{j}$, $j\in \NN_0$, and analogously for all other arising sequences. $\M$ is called {\itshape normalized} if $1=M_0\le M_1$ holds true.\vspace{6pt}

$\M$ is said to be {\itshape log-convex}, (for short, (lc)) if
$$\forall\;j\in\NN:\;M_j^2\le M_{j-1} M_{j+1},$$
equivalently if $\m$ is nondecreasing. If $\M$ is log-convex and normalized, then both $j\mapsto M_j$ and $j\mapsto(M_j)^{1/j}$ are nondecreasing and $(M_j)^{1/j}\le m_{j-1}$ for all $j\in\NN$. Finally $M_{j+k}\le M_j M_k$ follows for all $j,k\in\NN_0$.\vspace{6pt}

If $\check{\M}$ is log-convex, then $\M$ is called {\itshape strongly log-convex}, denoted by \hypertarget{slc}{$(\text{slc})$}. We say that a sequence $\M$ is a \emph{weight sequence} if it is (lc) and $\lim_{j\to\infty} m_j =\infty$.
We see that $\M$ is a \emph{normalized} weight sequence if and only if $1\leq m_0\le m_1\le\dots$, $\lim_{j\rightarrow+\infty}m_j=+\infty$ (e.g. see \cite[p. 104]{compositionpaper}) and there is a one-to-one correspondence between $\M$ and $\m$ by taking $M_j:=\prod_{i=0}^{j-1}m_i$.\vspace{6pt}

For $a\in\RR$ we set
$$\G^a:=(j!^a)_{j\in\NN_0},\hspace{15pt}\overline{\G}^a:=(j^{ja})_{j\in\NN_0},$$
i.e. for $a>0$ the sequence $\G^a$ is the Gevrey-sequence of index $a$. Clearly $ \G^a$ and $\overline{\G}^a $ are normalized weight sequences for any $a>0$ (by the convention $0^0:=1$).

%When $a=0$, then $G^0\equiv\overline{G}^0\equiv (1)_{j\in\NN_0}$.
%For any $M\in\hyperlink{LCset}{\mathcal{LC}}$ and $C\in\NN_{>0}$ we get that
%\begin{equation}\label{rootgeneralizedincr}
%\forall\;p\in\NN_0:\;\;\;(M_{Cp})^{1/C}\le(M_{(C+1)p})^{1/(C+1)},
%\end{equation}
%because this is equivalent to having $M_{Cp}(M_{Cp})^{1/C}\le M_{Cp+p}$ and so to $\mu_1\cdots\mu_{Cp}=M_{Cp}\le(\mu_{Cp+1}\cdots\mu_{Cp+p})^C$ which is satisfied since $p\mapsto\mu_p$ is nondecreasing.\vspace{6pt}

$\M$ satisfies the condition of {\itshape moderate growth}, denoted by \hypertarget{mg}{$(\text{mg})$}, if
$$\exists\;C\ge 1\;\forall\;j,k\in\NN_0:\;M_{j+k}\le C^{j+k} M_j M_k.$$
In the classical work of Komatsu~\cite{Komatsu73} this condition is named $(M.2)$ and also known in the literature under the name {\itshape stability under ultradifferential operators}. $\M$ satisfies the weaker requirement of {\itshape derivation closedness}, denoted by \hypertarget{dc}{$(\text{dc})$}, if
$$\exists\;D\ge 1\;\forall\;j\in\NN_0:\;M_{j+1}\le D^{j+1} M_j\Longleftrightarrow m_{j}\le D^{j+1}.$$
In \cite{Komatsu73} this is condition $(M.2')$. Both \hyperlink{mg}{$(\on{mg})$} and \hyperlink{dc}{$(\on{dc})$} are preserved when multiplying or dividing $\M$ by any sequence $\G^a$. In particular, both conditions hold simultaneously true or false for $\M$ and $\check{\M}$.

\vspace{6pt}

We say $\M$ {has the} {\itshape root almost increasing {property}}, denoted by \hypertarget{rai}{$(\on{rai})$}, if
the sequence of roots $({\check{M}}_j^{1/j})_{j\in\NN}$ is almost increasing, that is,
$$\exists\;C>0\;\forall\;1\le j\le k:\;\;\;\check{M}_j^{1/j}\le C\check{M}_k^{1/k}.$$
$\M$ has the {\itshape Fa\`{a}-di-Bruno property}, denoted by \hypertarget{FdB}{$(\on{FdB})$}, if
$$\exists\;C\ge 1\;\exists\;h\ge 1\;\forall\;j\in\NN_0:\;\;\;\check{M}^{\circ}_j\le Ch^j\check{M}_j,$$
where $\check{\M}^{\circ}:=(\check{M}_j^{\circ})_{j\in\NN_0}$ is the sequence defined by
\begin{equation}\label{circsequ}\check{M}_k^{\circ}:=\max\left\{\check{M}_\ell\cdot \check{M}_{j_1}\cdots \check{M}_{j_{\ell}}: j_i\in\NN, \sum_{i=1}^{\ell}j_i=k\right\},\ k\in\NN;\;\;\;\;\check{M}_0^{\circ}:=1.
\end{equation}
%$M$ is called {\itshape nonquasianalytic,} denoted by \hypertarget{mnq}{$(\text{nq})$}, if $M$ satisfies
%$$\sum_{j=1}^{+\infty}\frac{1}{\mu_j}<+\infty.$$

%We say that $M$ has condition \hypertarget{beta1}{$(\beta_1)$} (introduced in \cite{petzsche}), if
%$$\exists\;Q\in\NN_{\ge 2}:\;\;\;\liminf_{j\rightarrow+\infty}\frac{\mu_{Qj}}{\mu_j}>Q.$$
%In \cite{petzsche} the surjectivity of the {\itshape Borel mapping} $f\mapsto(f^{(j)}(0))_{j\in\NN_0}$ from $\mathcal{E}_{[M]}$ onto the corresponding weighted sequence space has been characterized in terms of this condition. More precisely, there for $M\in\hyperlink{LCset}{\mathcal{LC}}$ it has been shown that \hyperlink{beta1}{$(\beta_1)$} is equivalent to requiring \hyperlink{gamma1}{$(\gamma_1)$}, i.e.
%$$\sup_{j\in\NN_{>0}}\frac{\mu_j}{j}\sum_{k\ge j}\frac{1}{\mu_k}<+\infty.$$
%In the literature \hyperlink{gamma1}{$(\gamma_1)$} is also called ''strong nonquasianalyticity condition'' and in \cite{Komatsu73} it is denoted by $(M.3)$ (in fact, there $\frac{\mu_j}{j}$ is replaced by $\frac{\mu_j}{j-1}$ for $j\ge 2$ but which is equivalent to having \hyperlink{gamma1}{$(\gamma_1)$}).\vspace{6pt}

Let $\M,\LL\in\RR_{>0}^{\NN_0}$ be given {with arbitrary $M_0,L_0>0$}, we write $\M\hypertarget{preceq}{\preceq}\LL$ if $\sup_{j\in\NN}\left(M_j/L_j\right)^{1/j}<+\infty$ {or, equivalently, if there exist $A,B>0$ such that $M_j\le AB^jL_j$ for every $j\in\NN_0$}. We say $\M$ and $\LL$ are {\itshape equivalent}, denoted by $\M\hypertarget{approx}{\approx}\LL$, if $\M\hyperlink{preceq}{\preceq}\LL$ and $\LL\hyperlink{preceq}{\preceq}\M$.
Note that, in case $M_0=L_0=1$, equivalence amounts to $B^jM_j\le L_j\le C^jM_j$ for every $j\in\NN_0$ and suitable $B,C>0$.
Properties \hyperlink{mg}{$(\on{mg})$} and \hyperlink{dc}{$(\on{dc})$} are clearly preserved under \hyperlink{approx}{$\approx$}. %and for \hyperlink{beta1}{$(\beta_1)$} this follows by the characterizations obtained in \cite{petzsche}.

Let us write $\M\le \LL$ if $M_j\le L_j$ for all $j\in\NN_0$.\vspace{6pt}

Finally{,} we recall some useful elementary estimates,
%(a consequence of) Stirling's formula
\begin{equation}\label{Stirling}
	\forall\;j\in\NN:\;\;\;\frac{j^j}{e^j}\le j!\le j^j,
\end{equation}
which immediately imply that  %the following:
%
%\begin{lemma}\label{overlineGequiv}
%	We have
 $\G^a\hyperlink{approx}{\approx}\overline{\G}^a$ for any $a\in\RR$.
%\end{lemma}

\subsection{Associated weight function}\label{ssectAssocWeightFunct}
Let $\M\in\RR_{>0}^{\NN_0}$, then the {\itshape associated function} $\omega_\M: \RR_{\ge 0}\rightarrow\RR\cup\{+\infty\}$ is defined by
\begin{equation*}\label{assofunc}
\omega_\M(t):=\sup_{j\in\NN_0}\ln\left(\frac{t^j}{M_j}\right)\;\;\;\text{for}\;t>0,\hspace{30pt}\omega_\M(0):=0.
\end{equation*}
For an abstract introduction of the associated function we refer to \cite[Chapitre I]{mandelbrojtbook}, see also \cite[Definition 3.1]{Komatsu73}. %Note that $\omega_M$ is here extended to whole $\RR$ in an symmetric (even) way.

If $\liminf_{j\rightarrow+\infty}(M_j)^{1/j}>0$, then $\omega_\M(t)=0$ for sufficiently small $t>0$, since $t^0/M_0=1$ and $\ln\left(\frac{t^j}{M_j}\right)<0$ precisely if $t<(M_j)^{1/j}$, $j\in\NN$ (in particular, if $M_j\ge 1$ for all $j\in\NN_0$, then $\omega_\M$ vanishes on $[0,1]$). Moreover, under this assumption $t\mapsto\omega_\M(t)$ is a continuous nondecreasing function, which is convex in the variable $\ln(t)$ and tends faster to infinity than any $\ln(t^j)$, $j\ge 1$, as $t\rightarrow+\infty$. If $\lim_{j\rightarrow+\infty}(M_j)^{1/j}=+\infty$, then $\omega_\M(t)<+\infty$ for each finite $t$, so this will be a basic assumption for defining $\omega_\M$.\vspace{6pt}

If $\M$ is a weight sequence, then we can compute $\M$ by involving $\omega_\M$ as follows, see \cite[Chapitre I, 1.4, 1.8]{mandelbrojtbook} and also \cite[Prop. 3.2]{Komatsu73}:
\begin{equation}\label{Prop32Komatsu}
	M_j=\sup_{t\ge 0}\frac{t^j}{\exp(\omega_{\M}(t))},\;\;\;j\in\NN_0.
\end{equation}
Moreover, in this case one has
\begin{equation*}%\label{assovanishing}
	\omega_\M(t)=0\hspace{20pt}\forall\;t\in[0,m_0],
\end{equation*}
by the known integral representation formula for $\omega_\M$, see \cite[1.8. III]{mandelbrojtbook} and also \cite[$(3.11)$]{Komatsu73}.

If $\M\in\RR_{>0}^{\NN_0}$ satisfies $\lim_{j\rightarrow+\infty}(M_j)^{1/j}=+\infty$, then the right-hand side of formula \eqref{Prop32Komatsu} yields the $j$-th term of the log-convex minorant $\M^{\on{lc}}$ of $\M$, i.e. the log-convex sequence such that each log-convex sequence $\LL$ with $\LL\le \M$ satisfies $\LL\le \M^{\on{lc}}$ (moreover, $\M^{\on{lc}}\equiv \M$ if and only if $\M$ is log-convex). By the results from \cite[Chapitre I]{mandelbrojtbook} it also follows that $\omega_\M\equiv\omega_{\M^{\on{lc}}}$.

Finally, if for $\beta>0$ we write $\M^{1/\beta}:=(M_j^{1/\beta})_{j\in\NN_0}$, we recall the following immediate equality, e.g. see  \cite[$(2.7)$]{sectorialextensions}:
\begin{equation*}%\label{assofctpower}
	\forall\;t\ge 0:\;\;\;\omega_{\M}^{\beta}(t):=\omega_{\M}(t^{\beta})=\beta\omega_{\M^{1/\beta}}(t).
\end{equation*}

\subsection{Growth index $\gamma(\M)$}\label{ssect.IndexgammaM}

We say $\M$ satisfies property $\left(P_{\gamma}\right)$ if there exists a sequence of real numbers $\boldsymbol{\ell}=(\ell_{j})_{j\in\NN_0}$
such that:
\begin{enumerate}[(i)]
	\item $\m\simeq\boldsymbol{\ell}$, that is,
	$$\exists\;a\ge 1\;\forall\;j\in\NN_0:\;\;\;a^{-1}m_j\le \ell_j\le a m_j,$$
	
	\item $\left((j+1)^{-\gamma}\ell_{j}\right)_{j\in\NN_0}$ is nondecreasing.
\end{enumerate}
Note that $\m\simeq\boldsymbol{\ell}$ implies $\M\hyperlink{approx}{\approx}\LL$.

If $(P_\gamma)$ holds true for $\M$, then $(P_{\gamma'})$ also holds for any $\gamma'\leq\gamma$.
It is then natural to define the \textit{growth index} $\gamma(\M)$ by
$$\gamma(\M):=\sup\{\gamma\in\RR:\, (P_{\gamma})\hbox{ is fulfilled}\},$$%
with the conventions $\inf\emptyset=\sup\RR=+\infty$ and $\inf\RR=\sup\emptyset=-\infty$ (see \cite[Rem. 2.2]{index}). For a comprehensive study of this index we refer to \cite[Sect. 3]{index}, especially to the characterizing result \cite[Thm. 3.11]{index}.
This growth index was originally defined and considered for so-called {\itshape strongly regular sequences} by V. Thilliez in \cite[Sect. 1]{Thilliezdivision}.

\subsection{Weight functions}\label{ssect.WeightFunctions}
A function $\omega:[0,+\infty)\rightarrow[0,+\infty)$ is called a {\itshape weight function} (in the terminology of \cite[Sect. 2.1]{index}, \cite[Sect. 2.2]{sectorialextensions}, \cite[Sect. 2.2]{sectorialextensions1}), if it is continuous, nondecreasing, $\omega(0)=0$ and $\lim_{t\rightarrow+\infty}\omega(t)=+\infty$. If $\omega$ satisfies in addition $\omega(t)=0$ for all $t\in[0,1]$, then we call $\omega$ a {\itshape normalized weight function}. For convenience we will write that $\omega$ has $\hypertarget{om0}{(\omega_0)}$ if it is a normalized weight.\vspace{6pt}

For any $a>0$ we put $\omega^a$ for the function given by $\omega^a(t):=\omega(t^a)$, i.e. composing with a so-called Gevrey weight $t\mapsto t^a$.
%\red{(If $a=0$, then we put $\omega^0(t):=\omega(1)$.) Let us also set $\omega^{\iota}(t):=\omega(1/t)$ for $t>0$.}
\vspace{6pt}

%and finally, if $a<0$, then let us put $\omega^{a}:=(\omega^{\iota})^{-a}=(\omega^{-a})^{\iota}$ since in this case $\omega^a(t)=\omega(t^a)=\omega(1/t^{-a})=\omega^{\iota}(t^{-a})=(\omega^{\iota})^{-a}(t)=(\omega^{-a})^{\iota}(t)$ for all $t>0$.\vspace{6pt}

Let $\sigma,\tau$ be weight functions, we write $\sigma\hypertarget{ompreceq}{\preceq}\tau$ if $\tau(t)=O(\sigma(t))\;\text{as}\;t\rightarrow+\infty$
and call them equivalent, denoted by $\sigma\hypertarget{sim}{\sim}\tau$, if
$\sigma\hyperlink{ompreceq}{\preceq}\tau$ and $\tau\hyperlink{ompreceq}{\preceq}\sigma$.\vspace{6pt}

We consider the following (standard) conditions, this list of properties has already been used in~\cite{dissertation}.

\begin{itemize}
	\item[\hypertarget{om1}{$(\omega_1)}$] $\omega(2t)=O(\omega(t))$ as $t\rightarrow+\infty$, i.e. $\exists\;L\ge 1\;\forall\;t\ge 0:\;\;\;\omega(2t)\le L(\omega(t)+1)$.
	
	\item[\hypertarget{om2}{$(\omega_2)$}] $\omega(t)=O(t)$ as $t\rightarrow+\infty$.
	
	\item[\hypertarget{om3}{$(\omega_3)$}] $\ln(t)=o(\omega(t))$ as $t\rightarrow+\infty$.
	
	\item[\hypertarget{om4}{$(\omega_4)$}] $\varphi_{\omega}:t\mapsto\omega(e^t)$ is a convex function on $\RR$.
	
	\item[\hypertarget{om5}{$(\omega_5)$}] $\omega(t)=o(t)$ as $t\rightarrow+\infty$.
	
	\item[\hypertarget{om6}{$(\omega_6)$}] $\exists\;H\ge 1\;\forall\;t\ge 0:\;2\omega(t)\le\omega(H t)+H$.
\end{itemize}

%$\omega$ is called {\itshape non-quasianalytic}, denoted by \hypertarget{omnq}{$(\omega_{\on{nq}})$}, if
%$$\int_{1}^{+\infty}\frac{\omega(t)}{t^2}<+\infty,$$
%
%    NO NEED FOR CONDITION (SNQ)
%
%$\omega$ is called {\itshape strong non-quasianalytic}, denoted by \hypertarget{omsnq}{$(\omega_{\text{snq}})$}, if
%$$\exists\;C>0:\;\forall\;y>0: \int_1^{+\infty}\frac{\omega(y t)}{t^2}dt\le C\omega(y)+C.$$
For convenience we define the sets
$$\hypertarget{omset0}{\mathcal{W}_0}:=\{\omega:[0,\infty)\rightarrow[0,\infty): \omega\;\text{has}\;\hyperlink{om0}{(\omega_0)},\hyperlink{om3}{(\omega_3)},\hyperlink{om4}{(\omega_4)}\},\hspace{20pt}\hypertarget{omset1}{\mathcal{W}}:=\{\omega\in\mathcal{W}_0: \omega\;\text{has}\;\hyperlink{om1}{(\omega_1)}\}.$$
For any $\omega\in\hyperlink{omset0}{\mathcal{W}_0}$ we define the {\itshape Legendre-Fenchel-Young-conjugate} of $\varphi_{\omega}$ by
\begin{equation}\label{legendreconjugate}
	\varphi^{*}_{\omega}(x):=\sup\{x y-\varphi_{\omega}(y): y\ge 0\},\;\;\;x\ge 0,
\end{equation}
with the following properties, e.g. see \cite[Remark 1.3, Lemma 1.5]{BraunMeiseTaylor90}: It is convex and nondecreasing, $\varphi^{*}_{\omega}(0)=0$, $\varphi^{**}_{\omega}=\varphi_{\omega}$, $\lim_{x\rightarrow+\infty}\frac{x}{\varphi^{*}_{\omega}(x)}=0$ and finally $x\mapsto\frac{\varphi_{\omega}(x)}{x}$ and $x\mapsto\frac{\varphi^{*}_{\omega}(x)}{x}$ are nondecreasing on $[0,+\infty)$. Note that by normalization we can extend the supremum in \eqref{legendreconjugate} from $y\ge 0$ to $y\in\RR$ without changing the value of $\varphi^{*}_{\omega}(x)$ for given $x\ge 0$.\vspace{6pt}

Finally{,} let us introduce and recall the following crucial growth assumption on $\omega$:
\begin{equation}\label{alpha0}
	\exists\;C\ge 1\;\exists\;t_0\ge 0\;\forall\;\lambda\ge 1\;\forall\;t\ge t_0:\;\;\;\omega(\lambda t)\le C\lambda\omega(t).
\end{equation}
In the literature this condition is frequently denoted by $(\alpha_0)$. It is known that a weight function $\omega$ is equivalent to a subadditive weight function $\sigma$ (i.e., $\sigma(s+t)\le\sigma(s)+\sigma(t)$ for every $s,t\ge 0$), or even to a concave weight function, if and only if \eqref{alpha0} holds true, we refer to \cite[Sect. 4.1]{subaddlike} and the introduction of \cite{subaddlike} with the citations therein. In \cite[Thm. 4.5]{subaddlike} this condition for $\omega_{\M}$ has been characterized in terms of $\M$.

It is also known that $(\alpha_0)$ characterizes some desired stability properties for ultradifferentiable classes $\mathcal{E}_{[\omega]}$, e.g. closedness under composition, inverse closedness and closedness under solving ODE's. The definition of such classes (which will not be used in this paper) and these results can be found in \cite{compositionpaper}, \cite[Thm. 1, Thm. 3]{characterizationstabilitypaper} and \cite[Thm. 4.8]{almostanalytic} and in the references therein (see also~\cite{FernandezGalbis06} for closedness under composition).\vspace{6pt}

We recall the following known result, e.g. see \cite[Sect. 5]{compositionpaper} and \cite[Lemma 4.1]{Komatsu73}, \cite[Lemma 2.8]{testfunctioncharacterization} and \cite[Lemma 2.4]{sectorialextensions} and the references mentioned in the proofs there.

\begin{lemma}\label{assoweightomega0}
	Let $\M$ be a normalized weight sequence, then $\omega_\M\in\hyperlink{omset0}{\mathcal{W}_0}$ holds true. Moreover,
	\begin{itemize}
		\item[$(i)$] $\liminf_{j\rightarrow\infty}(\check{M}_j)^{1/j}>0$ if and only if \hyperlink{om2}{$(\omega_2)$} holds for $\omega_\M$,
		
		\item[$(ii)$] $\lim_{j\rightarrow\infty}(\check{M}_j)^{1/j}=+\infty$ if and only if \hyperlink{om5}{$(\omega_5)$} holds for $\omega_\M$,
		
		\item[$(iii)$] \hyperlink{om6}{$(\omega_6)$} holds for $\omega_\M$ if and only if $\M$ does have \hyperlink{mg}{$(\on{mg})$}.
		
		%\item[$(iv)$] \hyperlink{omnq}{$(\omega_{\on{nq}})$} for $\omega_M$ if and only if $M$ does have \hyperlink{mnq}{$(\on{nq})$}.
	\end{itemize}
\end{lemma}

\subsection{Weight matrices}\label{ssect.WeightMatrices}
For the following definitions and conditions see also \cite[Sect. 4]{compositionpaper}.

Let $\mathcal{I}=\RR_{>0}$ denote the index set (equipped with the natural order), a {\itshape weight matrix} $\mathcal{M}$ associated with $\mathcal{I}$ is a (one parameter) family of sequences $\mathcal{M}:=\{\M^{(\alpha)}\in\RR_{>0}^{\NN_0}: \alpha\in\mathcal{I}\}$, such that
$$\M^{(\alpha)}\le \M^{(\beta)}\;\text{for}\;\alpha\le\beta;\hspace{15pt}M^{(\alpha)}_0=1,\;\;\;\forall\;\alpha\in\mathcal{I}.$$
We call a weight matrix $\mathcal{M}$ {\itshape log-convex,} denoted by \hypertarget{Mlc}{$(\mathcal{M}_{\on{lc}})$}, if $\M^{(\alpha)}$ is a log-convex sequence for all $\;\alpha\in\mathcal{I}$.
Moreover, we say that a weight matrix $\mathcal{M}$ is {\itshape standard log-convex,} abbreviated by \hypertarget{Msc}{$(\mathcal{M}_{\on{sc}})$}, if $\M^{(\alpha)}$ is a normalized weight sequence for all $\;\alpha\in\mathcal{I}$.
%$$\forall\;\alpha\in\mathcal{I}:\;\M^{(\alpha)}\in\hyperlink{LCset}{\mathcal{LC}}.$$
We put $\check{M}^{(\alpha)}_j:=\frac{M^{(\alpha)}_j}{j!}$ for $j\in\NN_0$, and $m^{(\alpha)}_j:=\frac{M^{(\alpha)}_{j+1}}{M^{(\alpha)}_{j}}$ for $j\in\NN_0$.

If $\mathcal{M}$ is a weight matrix with $\lim_{j\rightarrow\infty}(M^{(\alpha)}_j)^{1/j}=+\infty$ for all $\alpha$, then let us set
\begin{equation*}%\label{lcmatrix}
	\mathcal{M}^{\on{lc}}:=\{(\M^{(\alpha)})^{\on{lc}}: \M^{(\alpha)}\in\mathcal{M}\}.
\end{equation*}
For $\alpha\le\beta$, since $\M^{(\alpha)}\le \M^{(\beta)}$ we have $(\M^{(\alpha)})^{\on{lc}}\le(\M^{(\beta)})^{\on{lc}}$. Moreover, $(M^{(\alpha)})^{\on{lc}}_0=M^{(\alpha)}_0=1$.

A matrix is called {\itshape constant} if $\M^{(\alpha)}\hyperlink{approx}{\approx}\M^{(\beta)}$ for all $\alpha,\beta\in\mathcal{I}$.\vspace{6pt}

Let $\mathcal{M}=\{\M^{(\alpha)}: \alpha\in\mathcal{I}\}$ and $\mathcal{L}=\{\LL^{(\alpha)}: \alpha\in\mathcal{I}\}$ be given. We write $\mathcal{M}\hypertarget{Mroumpreceq}{\{\preceq\}}\mathcal{L}$ if
$$\forall\;\alpha\in\mathcal{I}\;\exists\;\beta\in\mathcal{I}:\;\;\;\M^{(\alpha)}\hyperlink{preceq}{\preceq}\LL^{(\beta)},$$
and call $\mathcal{M}$ and $\mathcal{L}$ $R$-equivalent, if
$\mathcal{M}\hyperlink{Mroumpreceq}{\{\preceq\}}\mathcal{L}$ and $\mathcal{L}\hyperlink{Mroumpreceq}{\{\preceq\}}\mathcal{M}$.\vspace{6pt}

%A matrix is called {\itshape non-quasianalytic} if any sequence $M^{(\alpha)}$ is non-quasianalytic. When dealing with Roumieu type classes then it suffices to assume that there exists $\alpha_0\in\mathcal{I}$ such that $M^{(\alpha_0)}$ is non-quasianalytic since smaller indices can be skipped, see also the discussion in \cite[Sect. 5.1]{borelmappingquasianalytic}.\vspace{6pt}

Let us consider the following crucial assumptions (of Roumieu-type) on a given weight matrix $\mathcal{M}$, see \cite[Sect. 4.1]{compositionpaper} and \cite[Sect. 7.2]{dissertation}:\vspace{6pt}

\hypertarget{R-Comega}{$(\mathcal{M}_{\{\text{C}^{\omega}\}})$} \hskip1cm $\exists\;\alpha\in\mathcal{I}:\;\;\;\liminf_{j\rightarrow\infty}(\check{M}^{(\alpha)}_j)^{1/j}>0$,

\hypertarget{holom}{$(\mathcal{M}_{\mathcal{H}})$} \hskip1cm $\forall\;\alpha\in\mathcal{I}:\;\;\;\liminf_{j\rightarrow\infty}(\check{M}_j^{(\alpha)})^{1/j}>0$,

\hypertarget{R-rai}{$(\mathcal{M}_{\{\text{rai}\}})$} \hskip1cm $\forall\;\alpha\in\mathcal{I}\;\exists\;C>0\;\exists\;\beta\in\mathcal{I}\;\forall\;1\le j\le k:\;\;\;(\check{M}^{(\alpha)}_j)^{1/j}\le C(\check{M}^{(\beta)}_k)^{1/k}$,

\hypertarget{R-FdB}{$(\mathcal{M}_{\{\text{FdB}\}})$} \hskip1cm  $\forall\;\alpha\in\mathcal{I}\;\exists\;\beta\in\mathcal{I}:\;\;\; (\check{\M}^{(\alpha)})^{\circ}\hyperlink{mpreceq}{\preceq} \check{\M}^{(\beta)}$,

where $(\check{\M}^{(\alpha)})^{\circ}$
%:=((\check{M}^{(\alpha)}_j)^{\circ})_j$
is the sequence defined by \eqref{circsequ}.
%$$(m^{(\alpha)}_k)^{\circ}:=\max\left\{m^{(\alpha)}_j\cdot m^{(\alpha)}_{j_1}\cdots m^{(\alpha)}_{j_{\ell}}: j_i\in\NN_{>0}, \sum_{i=1}^{\ell}j_i=k\right\},\;\;\;\;\;(m^{(\alpha)}_0)^{\circ}:=1.$$

Moreover, let us consider

\hypertarget{R-mg}{$(\mathcal{M}_{\{\text{mg}\}})$} \hskip1cm $\forall\;\alpha\in\mathcal{I}\;\exists\;C>0\;\exists\;\beta\in\mathcal{I}\;\forall\;j,k\in\NN_0: M^{(\alpha)}_{j+k}\le C^{j+k} M^{(\beta)}_j M^{(\beta)}_k$,\par\vskip.3cm
and the weaker requirement

\hypertarget{R-dc}{$(\mathcal{M}_{\{\text{dc}\}})$} \hskip1cm $\forall\;\alpha\in\mathcal{I}\;\exists\;C>0\;\exists\;\beta\in\mathcal{I}\;\forall\;j\in\NN_0: M^{(\alpha)}_{j+1}\le C^{j+1} M^{(\beta)}_j$.\par\vskip.3cm

%\subsection{Preliminary results}
Let us gather now some relevant information needed in the forthcoming sections.

\begin{lemma} Let $\mathcal{M}=\{\M^{(\alpha)}: \alpha\in\mathcal{I}\}$ be a weight matrix. If $\mathcal{M}$ has \hyperlink{R-rai}{$(\mathcal{M}_{\{\on{rai}\}})$}, then
	\begin{multline}
	    \label{theorem1siddiqiequ0}
		\forall\;\alpha>0\;\exists\;H\ge 1\;\exists\;\alpha'(\ge \alpha)\;\forall\;k\in\NN\;\forall\;j_1,\dots,j_k\in\NN_0:\\\check{M}^{(\alpha)}_{j_1}\cdots \check{M}^{(\alpha)}_{j_k}\le H^{j_1+\dots+j_k} \check{M}^{(\alpha')}_{j_1+\dots+j_k}.
	\end{multline}
Note that the indices $\alpha$ and $\alpha'$ are related by property \hyperlink{R-rai}{$(\mathcal{M}_{\{\on{rai}\}})$}.
\end{lemma}
\demo{Proof}
If $j_1,\dots,j_k\ge 1$ we estimate by
\begin{align*}
    \check{M}^{(\alpha)}_{j_1}\cdots \check{M}^{(\alpha)}_{j_k}&\le H^{j_1}(\check{M}^{(\alpha')}_{j_1+\dots+j_k})^{\frac{j_1}{j_1+\dots +j_k}}\cdots H^{j_k}(\check{M}^{(\alpha')}_{j_1+\dots+j_k})^{\frac{j_k}{j_1+\dots +j_k}}\\
&=H^{j_1+\dots+j_k}\check{M}^{(\alpha')}_{j_1+\dots+j_k},
\end{align*}
and the remaining cases follow by $\check{M}^{(\alpha)}_0=M^{(\alpha)}_0=1$.
\qed\enddemo

\begin{lemma}\label{condcomparison}
	Let $\mathcal{M}=\{\M^{(\alpha)}: \alpha\in\mathcal{I}\}$ be a weight matrix. Then we have the following:
	\begin{itemize}
		\item[$(i)$] \hyperlink{R-rai}{$(\mathcal{M}_{\{\on{rai}\}})$} implies \hyperlink{holom}{$(\mathcal{M}_{\mathcal{H}})$}.
		
		\item[$(ii)$] \hyperlink{R-dc}{$(\mathcal{M}_{\{\on{dc}\}})$} and \hyperlink{R-rai}{$(\mathcal{M}_{\{\on{rai}\}})$} imply \hyperlink{R-FdB}{$(\mathcal{M}_{\{\on{FdB}\}})$}.
		
		\item[$(iii)$] If
		\begin{equation}\label{rai}
			\forall\;\alpha\in\mathcal{I}\;\exists\;H\ge 1\;\forall\;1\le j\le k:\;\;\;(M^{(\alpha)}_j)^{1/j}\le H(M^{(\alpha)}_k)^{1/k},
		\end{equation}
		i.e. each sequence $((M^{(\alpha)}_j)^{1/j})_j$ is almost increasing, then \hyperlink{holom}{$(\mathcal{M}_{\mathcal{H}})$} and \hyperlink{R-FdB}{$(\mathcal{M}_{\{\on{FdB}\}})$} imply \hyperlink{R-rai}{$(\mathcal{M}_{\{\on{rai}\}})$}.
		
		In particular, \eqref{rai} holds true (with $H=1$ for any $\alpha$) provided that $\mathcal{M}$ is log-convex.
	\end{itemize}
\end{lemma}

\demo{Proof}
$(i)$ By the order of the sequences we can assume w.l.o.g. $\beta\ge\alpha$ and for each $\alpha\in\mathcal{I}$ there exists a minimal $\beta=\beta(\alpha)\ge\alpha$ such that $\check{\M}^{(\alpha)}$ and $\check{\M}^{(\beta)}$ are related by \hyperlink{R-rai}{$(\mathcal{M}_{\{\on{rai}\}})$}. Then $(\check{M}^{(\beta)}_j)^{1/j}\ge\frac{\check{M}^{(\alpha)}_1}{C}>0$ for some $C\ge 1$ and all $j\ge 1$ (see also \cite[Lemma 3.6 $(ii)$]{subaddlike}). Since w.l.o.g. we can restrict in the Roumieu case to all $\beta(\alpha)$ (yielding an $R$-equivalent matrix) we are done.\vspace{6pt}

$(ii)$ See the proofs of \cite[Thm. 4.9 $(3)\Rightarrow(4)$]{compositionpaper} and \cite[Lemma 8.2.3 $(2)$]{dissertation}.

$(iii)$ See the proofs of \cite[Lemma 1 $(2)$]{characterizationstabilitypaper} and \cite[Lemma 8.2.3 $(4)$]{dissertation}.
\qed\enddemo

\subsection{Weight matrices associated with weight functions}\label{ssect.WeightMatricesFromFunctions}
We summarize some facts which are shown in \cite[Section 5]{compositionpaper} and are needed in this work. All properties listed below are valid for $\omega\in\hyperlink{omset0}{\mathcal{W}_0}$, except \eqref{newexpabsorb} for which \hyperlink{om1}{$(\omega_1)$} is necessary.

\begin{itemize}
	\item[$(i)$] The idea was that to each $\omega\in\hyperlink{omset0}{\mathcal{W}_0}$ we can associate a standard log-convex weight matrix $\mathcal{M}_{\omega}:=\{\W^{(\ell)}=(W^{(\ell)}_j)_{j\in\NN_0}: \ell>0\}$ by\vspace{6pt}
	
	 \centerline{$W^{(\ell)}_j:=\exp\left(\frac{1}{\ell}\varphi^{*}_{\omega}(\ell j)\right)$.}\vspace{6pt}
	
	%The corresponding sequences of quotients are denoted by $\vartheta^{(l)}$ and we even have $\vartheta^{(l_1)}\le\vartheta^{(l_2)}$ for $l_1\le l_2$, see \cite[Sect. 2.5]{whitneyextensionweightmatrix}
	%In general it is not clear that $W^x$ is strongly log-convex, i.e. $w^x$ is log-convex, too.
	
	\item[$(ii)$] $\mathcal{M}_{\omega}$ satisfies
	\begin{equation}\label{newmoderategrowth}
		 \forall\;\ell>0\;\forall\;j,k\in\NN_0:\;\;\;W^{(\ell)}_{j+k}\le W^{(2\ell)}_jW^{(2\ell)}_k,
	\end{equation}
	so both \hyperlink{R-mg}{$(\mathcal{M}_{\{\on{mg}\}})$} and \hyperlink{R-dc}{$(\mathcal{M}_{\{\on{dc}\}})$} are satisfied.
	
	\item[$(iii)$] \hyperlink{om6}{$(\omega_6)$} holds if and only if some/each $\W^{(\ell)}$ satisfies \hyperlink{mg}{$(\on{mg})$} if and only if $\W^{(\ell)}\hyperlink{approx}{\approx}\W^{(\ell_1)}$ for each $\ell,\ell_1>0$. Consequently \hyperlink{om6}{$(\omega_6)$} characterizes the situation when $\mathcal{M}_{\omega}$ is constant.
	
	\item[$(iv)$] In case $\omega$ has in addition \hyperlink{om1}{$(\omega_1)$}, then $\mathcal{M}_{\omega}$ has also
	\begin{equation}\label{newexpabsorb}
		\forall\;h\ge 1\;\exists\;A\ge 1\;\forall\;\ell>0\;\exists\;D\ge 1\;\forall\;j\in\NN_0:\;\;\;h^jW^{(\ell)}_j\le D W^{(A\ell)}_j,
	\end{equation}
see \cite[Lemma 5.9 (5.10)]{compositionpaper}.
%	and this estimate is crucial for proving $\mathcal{E}_{\{\mathcal{M}_{\omega}\}}=\mathcal{E}_{\{\omega\}}$ as locally convex vector spaces.
	 %$\mathcal{E}_{[\mathcal{M}_{\omega}]}=\mathcal{E}_{[\omega]}$ as locally convex vector spaces.
	
	\item[$(v)$] We have $\omega\hyperlink{sim}{\sim}\omega_{\W^{(\ell)}}$ for each $\ell>0$, more precisely
	\begin{equation}\label{goodequivalenceclassic}
		\forall\;\ell>0\,\,\exists\,D_{\ell}>0\;\forall\;t\ge 0:\;\;\;\ell\omega_{\W^{(\ell)}}(t)\le\omega(t)\le 2\ell\omega_{\W^{(\ell)}}(t)+D_{\ell},
	\end{equation}
	see \cite[Thm. 4.0.3, Lemma 5.1.3]{dissertation} and also \cite[Lemma 2.5]{sectorialextensions}.
	
	\item[$(vi)$] $\mathcal{M}_{\omega}$ satisfies \hyperlink{holom}{$(\mathcal{M}_{\mathcal{H}})$} if and only if $\omega$ has in addition \hyperlink{om2}{$(\omega_2)$} (by \eqref{goodequivalenceclassic} and $(i)$ in Lemma \ref{assoweightomega0}).
	
	\item[$(vii)$] Lemma \ref{condcomparison} (\cite[Lemma 1]{characterizationstabilitypaper}) applies to $\mathcal{M}_{\omega}$; so in view of $(i)$, $(ii)$ and $(vi)$ we have that for $\omega\in\hyperlink{omset0}{\mathcal{W}_0}$ with \hyperlink{om2}{$(\omega_2)$} properties \hyperlink{R-rai}{$(\mathcal{M}_{\{\on{rai}\}})$} and \hyperlink{R-FdB}{$(\mathcal{M}_{\{\on{FdB}\}})$} for $\mathcal{M}_{\omega}$ are simultaneously satisfied or violated.
\end{itemize}

\subsection{The growth index $\gamma(\omega)$}\label{ssect.IndexGammaOmega}

Let $\omega$ be a weight function. We recall the definition of the growth index $\gamma(\omega)$, see \cite[Sect. 2.3]{index} and the references therein: Let $\gamma>0$, then we say that $\omega$ has property $(P_{\omega,\gamma})$ if
\begin{equation*}\label{newindex1}
	 \exists\;K>1:\;\;\;\limsup_{t\rightarrow+\infty}\frac{\omega(K^{\gamma}t)}{\omega(t)}<K.
\end{equation*}
If $(P_{\omega,\gamma})$ holds for some $K>1$, then also $(P_{\omega,\gamma'})$ is satisfied for all $\gamma'\le\gamma$ with the same $K$. Moreover we can restrict to $\gamma>0$, because for $\gamma\le 0$ condition $(P_{\omega,\gamma})$ is satisfied for all weights $\omega$ (since $\omega$ is nondecreasing and $K>1$). Then we put
\begin{equation*}%\label{newindex2}
	\gamma(\omega):=\sup\{\gamma>0: (P_{\omega,\gamma})\;\;\text{is satisfied}\}.
\end{equation*}

We recall some facts about $\gamma(\omega)$:

\begin{itemize}
	\item[$(i)$] If $\omega\hyperlink{sim}{\sim}\sigma$ then $\gamma(\omega)=\gamma(\sigma)$, see \cite[Rem. 2.12]{index}.
	
	\item[$(ii)$] $\gamma(\omega)>0$ holds if and only if \hyperlink{om1}{$(\omega_1)$}, see \cite[Cor. 2.14]{index}.
	
%	\red{\item[$(-)$] $\gamma(\omega)>1$ holds if and only if \hyperlink{omsnq}{$(\omega_{\text{snq}})$}, see \cite[Cor. 2.13]{index}.}
	
	\item[$(iii)$] By definition one has $\gamma(\omega^a)=\frac{1}{a}\gamma(\omega)$ for any $a>0$.
	
	\item[$(iv)$] If $\omega\in\hyperlink{omset0}{\mathcal{W}_0}$ is given with associated weight matrix $\mathcal{M}_{\omega}:=\{\W^{(\ell)}: \ell>0\}$ and $\gamma(\omega)>\beta$, then \eqref{goodequivalenceclassic} implies $\gamma(\omega_{\W^{(\ell)}})>\beta$, but in general only $\gamma(\W^{(\ell)})\le\gamma(\omega_{\W^{(\ell)}})$ by \cite[Cor. 4.6 $(i)$]{index}. Here $\gamma(\W^{(\ell)})$ is the index in Subsection~\ref{ssect.IndexgammaM}, see more details in \cite[Sect. 3]{index}.
\end{itemize}

\section{Ultraholomorphic classes}\label{sect.UltraholomClass}
We introduce now the classes under consideration in this paper, see also \cite[Sect. 2.5]{sectorialextensions} and \cite[Sect. 2.5]{sectorialextensions1}. For the following definitions, notation and more details we refer to \cite[Section 2]{Sanzflatultraholomorphic}. Let $\mathcal{R}$ be the Riemann surface of the logarithm. We wish to work in general unbounded sectors in $\mathcal{R}$ with vertex at $0$, but all our results will be unchanged under rotation, so we will only consider sectors bisected by direction $0$: For $\alpha>0$ we set \[S_{\alpha}:=\left\{z\in\mathcal{R}: |\arg(z)|<\frac{\alpha\pi}{2}\right\},\]
i.e. the unbounded sector of opening $\alpha\pi$, bisected by direction $0$.

Let $\M$ be a sequence, $S\subseteq\mathcal{R}$ an (unbounded) sector and $h>0$. We define
\[\mathcal{A}_{\M,h}(S):=\{f\in\mathcal{H}(S): \|f\|_{\M,h}:=\sup_{z\in S, j\in\NN_0}\frac{|f^{(j)}(z)|}{h^j M_j}<+\infty\}.\]
$(\mathcal{A}_{\M,h}(S),\|\cdot\|_{\M,h})$ is a Banach space and we put
\begin{equation*}%\label{spacedefine}
	\mathcal{A}_{\{\M\}}(S):=\bigcup_{h>0}\mathcal{A}_{\M,h}(S).
\end{equation*}
$\mathcal{A}_{\{\M\}}(S)$ is called the Denjoy-Carleman ultraholomorphic class (of Roumieu type) associated with $\M$ in the sector $S$ (it is an $(LB)$ space). By definition it is immediate that $\M\hyperlink{approx}{\approx}\LL$ implies $\mathcal{A}_{\{\M\}}(S)=\mathcal{A}_{\{\LL\}}(S)$ (as locally convex vector spaces) for any sector $S$.
%
%
%      we supress the classes with uniform asymptotic expansion
%    THIS SHOULD BE THE TOPIC OF ANOTHER PAPER
%
%

%We also consider the set
%\[\mathcal{\widetilde{A}}^{u}_{\M,h}(S):=\{f\in\mathcal{H}(S): \left\|f\right\|_{\M,h,\overset{\sim}{u}}:=\sup_{z\in S,j\in\NN_{0}}\frac{|f(z)-\sum_{k=0}^{j-1}a_kz^k|}{h^{j}M_{j}|z|^j}<+\infty\}.\]
%$(\mathcal{\widetilde{A}}^{u}_{\M,h}(S),\|\cdot\|_{\M,h,\overset{\sim}{u}})$ is a Banach space and we put
%\begin{equation*}%\label{spacedefine}
%		 \mathcal{\widetilde{A}}^{u}_{\{\M\}}(S):=\bigcup_{h>0}\mathcal{\widetilde{A}}^{u}_{\M,h}(S).
%\end{equation*}
%$\widetilde{\mathcal{A}}^u_{\{\M\}}(S)$ stands for the $(LB)$ space of functions admitting a uniform $\{\M\}$-asymptotic expansion in $S$.
%By definition it is immediate that $\M\hyperlink{approx}{\approx}\LL$ implies $\mathcal{A}_{\{\M\}}(S)=\mathcal{A}_{\{\LL\}}(S)$ and $\mathcal{\widetilde{A}}^{u}_{\{\M\}}(S)=\mathcal{\widetilde{A}}^{u}_{\{\LL\}}(S)$ (as locally convex vector spaces) for any sector $S$.
\vspace{6pt}

Similarly as for the ultradifferentiable case, we now define ultraholomorphic classes associated with $\omega\in\hyperlink{omset0}{\mathcal{W}_0}$. Given an unbounded sector $S$, and for every $\ell>0$, we first define
\[\mathcal{A}_{\omega,\ell}(S):=\{f\in\mathcal{H}(S): \|f\|_{\omega,\ell}:=\sup_{z\in S, j\in\NN_0}\frac{|f^{(j)}(z)|}{\exp(\frac{1}{\ell}\varphi^{*}_{\omega}(\ell j))}<+\infty\}.\]
$(\mathcal{A}_{\omega,\ell}(S),\|\cdot\|_{\omega,\ell})$ is a Banach space and we put
\begin{equation*}%\label{spacedefine}
	 \mathcal{A}_{\{\omega\}}(S):=\bigcup_{\ell>0}\mathcal{A}_{\omega,\ell}(S).
\end{equation*}
$\mathcal{A}_{\{\omega\}}(S)$ is called the Denjoy-Carleman ultraholomorphic class (of Roumieu type) associated with $\omega$ in the sector $S$ (it is an $(LB)$ space). Again, equivalent weight functions provide equal associated ultraholomorphic classes.

Finally, we define ultraholomorphic classes of {\itshape Roumieu type} defined by a weight matrix $\mathcal{M}$ analogously as the ultradifferentiable counterparts introduced in \cite[Section 7]{dissertation} and also in \cite[Section 4.2]{compositionpaper}.

Given a weight matrix $\mathcal{M}=\{\M^{(\alpha)}\in\RR_{>0}^{\NN_0}: \alpha\in\mathcal{I}\}$ and a sector $S$ we may introduce the class $\mathcal{A}_{\{\mathcal{M}\}}(S)$  of {\itshape Roumieu type} as
\begin{equation*}%\label{generalroumieuultraholom}
	 \mathcal{A}_{\{\mathcal{M}\}}(S):=\bigcup_{\alpha\in\mathcal{I}}\mathcal{A}_{\{\M^{(\alpha)}\}}(S).
\end{equation*}
R-equivalent weight matrices yield (as locally convex vector spaces) the same function class on each sector $S$.

Let now $\omega\in\hyperlink{omset1}{\mathcal{W}}$ be given and let $\mathcal{M}_{\omega}$ be the associated weight matrix defined in Subsection \ref{ssect.WeightMatricesFromFunctions}, then
\begin{equation}\label{equaEqualitySpacesWeightFunctionMatrix}
	 \mathcal{A}_{\{\omega\}}(S)=\mathcal{A}_{\{\mathcal{M}_{\omega}\}}(S)
\end{equation}
holds as locally convex vector spaces. This equality is an easy consequence of \eqref{newexpabsorb} and the way the seminorms are defined in these spaces.\vspace{6pt}

On the other hand, by $(iii)$ in Subsection \ref{ssect.WeightMatricesFromFunctions} we get the following:

\begin{lemma}\label{om6lemma}
	Let $\omega\in\hyperlink{omset1}{\mathcal{W}}$ be given and assume that $\omega$ has \hyperlink{om6}{$(\omega_6)$}. Then, for all sectors $S$ we get that
$$\forall\;\ell>0:\;\;\;\mathcal{A}_{\{\omega\}}(S)=\mathcal{A}_{\{\W^{(\ell)}\}}(S)$$
as locally convex vector spaces.
\end{lemma}

If $f$ belongs to any of such classes, we may define the complex numbers $$f^{(j)}(0):=\lim_{z\in S, z\rightarrow 0}f^{(j)}(z),\quad j\in\NN_0.
$$

\section{Characteristic functions in ultraholomorphic classes}\label{sect.CharactFunctions}
%\subsection{Preliminaries}

We start with the following definition.

%
%
%    IN THIS SECTION WE HAVE SUPRESSED THE INFORMATION FOR CLASSES
%     WITH UNIFORM ASYMPTOTIC EXPANSION
%
%

\begin{definition}\label{chardef}
	Let $\LL\in\RR_{>0}^{\NN_0}$ and $S$ be a given sector.
%A function $f\in\widetilde{\mathcal{A}}^u_{\{\LL\}} (S)$ is said to be \textbf{characteristic} in the class $\widetilde{\mathcal{A}}^u_{\{\LL\}} (S)$, if the following holds true: If for some $\M\in\RR_{>0}^{\NN_0}$ we have $f\in\widetilde{\mathcal{A}}^u_{\{\M\}} (S)\subseteq\widetilde{\mathcal{A}}^u_{\{\LL\}} (S)$, then already $\widetilde{\mathcal{A}}^u_{\{\M\}} (S)= \widetilde{\mathcal{A}}^u_{\{\LL\}} (S)$.
%	
%	Similarly,
A function $f\in \mathcal{A}_{\{\LL\}} (S)$ is said to be \textit{characteristic} in the class $\mathcal{A}_{\{\LL\}} (S)$ if, whenever $f\in \mathcal{A}_{\{\M\}} (S)\subseteq \mathcal{A}_{\{\LL\}} (S)$ for some $\M\in\RR_{>0}^{\NN_0}$, we have that $\mathcal{A}_{\{\M\}} (S)= \mathcal{A}_{\{\LL\}} (S)$.
\end{definition}

%Let us fix some notation. If $f\in \widetilde{\mathcal{A}}^u_{\{\LL\}} (S)$  with $f\sim^u_{\{\LL\}} \sum_{j=0}^\infty c_j z^j$, then put
%$$r_f(z,n):= \frac{f(z)-\sum_{j=0}^{n-1} c_j z^j}{z^n},\;\;\;z\in S,\;n\in\NN_0.$$
%We observe that $z r_f(z,n+1)=r_f(z,n)-c_n$. If $f\in\widetilde{\mathcal{A}}^u_{\{\LL\}} (S)$, then $z\mapsto r_f(z,n)$ is bounded in $S$ for each $n\in\NN_0$ fixed and so
%$$\forall\;n\in\NN_0:\;\;\;\lim_{z\in S,z\rightarrow 0} r_f(z,n)=c_n.$$
%We consider the sequence defined by
%\begin{equation}\label{tildesequenceC}
%	\widetilde{C}_n(f):=\sup_{z\in S}|r_f(z,n)|, \qquad n\in\NN_0.
%\end{equation}
%
%\begin{theorem}\label{optthm1}
%	Let $\LL\in\RR_{>0}^{\NN_0}$ and $S$ be a given sector. Let $f\in\widetilde{\mathcal{A}}^u_{\{\LL\}} (S)$ with $f\sim^u_{\{\LL\}} \sum_{j=0}^\infty c_j z^j$, then each condition implies the following one:
%	\begin{enumerate}
%		\item The sequence $(|c_j|)_{j\in\NN_0}$ is equivalent to $\LL$.
%		\item The sequence $(|\widetilde{C}_j(f)|)_{j\in\NN_0}$ is equivalent to $\LL$.
%		\item $f$ is characteristic in the class $\widetilde{\mathcal{A}}^u_{\{\LL\}} (S)$.
%	\end{enumerate}
%\end{theorem}
%
%\demo{Proof}
%(1) $\Rightarrow$ (2) The estimate from above for $\widetilde{C}_n(f)$ follows by definition and the estimate from below holds since $\widetilde{C}_n(f)\geq |c_n|$.
%
%(2)  $\Rightarrow$ (3) This is clear.
%\qed\enddemo

For $f\in \mathcal{A}_{\{\LL\}}(S)$ we consider the sequence defined by
\begin{equation*}%\label{sequenceC}
	C_n(f):=\sup_{z\in S} |f^{(n)}(z)|, \qquad n\in\NN_0.
\end{equation*}
The next statement provides conditions on $f$ which imply it is characteristic.
%Similarly as before we get:

\begin{theorem}\label{optthm2}
	Let $\LL\in\RR_{>0}^{\NN_0}$, $S$ be a given sector and $f\in \mathcal{A}_{\{\LL\}} (S)$. Then, each of the following  conditions implies the next one:
	\begin{enumerate}
		\item[$(1)$] The sequence $(|f^{(j)}(0)|)_{j\in\NN_0}$ is equivalent to $\LL$.
		\item[$(2)$] The sequence $(C_j(f))_{j\in\NN_0}$ is equivalent to $\LL$.
		\item[$(3)$] $f$ is characteristic in the class $\mathcal{A}_{\{\LL\}} (S)$.
	\end{enumerate}
\end{theorem}

\demo{Proof}
(1) $\Rightarrow$ (2) As $f\in \mathcal{A}_{\{\LL\}} (S)$, there exist $A,B>0$ such that $C_n(f)\le AB^nL_n$ for every $n\in\NN_0$. On the other hand, it is clear that $C_n(f)\geq |f^{(n)}(0)|$, and the hypothesis allows us to conclude the other estimate.

(2)  $\Rightarrow$ (3) By assumption, there exist $A,B>0$ such that $L_n\le A{B}^nC_n(f)$ for every $n\in\NN_0$. If for some $\M=(M_n)_{n\in\NN_0}\in\RR_{>0}^{\NN_0}$ we have $f\in \mathcal{A}_{\{\M\}} (S)\subseteq \mathcal{A}_{\{\LL\}} (S)$, there exist $C,D>0$ such that $C_n(f)\le C {D}^nM_n$ for every $n\in\NN_0$. The two deduced inequalities show that $L_n\le{ AC(BD)^n}M_n$ for every $n\in\NN_0$, what easily implies that $\mathcal{A}_{\{\LL\}} (S)\subseteq \mathcal{A}_{\{\M\}} (S)$, and we are done.
\qed\enddemo

\subsection{Basic functions}\label{ssect.BasicFct}

Recall the notations {$\G^s:=(j!^s)_{j\in\NN_0}$ and }$\overline{\G}^s:=(j^{js})_{j\in\NN_0}$, $s\in\RR$, and that $\overline{\G}^s\hyperlink{approx}{\approx}\G^s$, see \eqref{Stirling}.\vspace{6pt}

%In fact, we need in the following $s\in\RR$ (so also negative $s$ are allowed) but

The \textit{{two-parametric} Mittag-Leffler function} is defined for all {complex parameters $A, B$ with $\Re(A)>0$} by
$$E_{A,B}(z):=\sum_{j=0}^\infty \frac{z^j}{\Gamma (Aj+B)},\quad z\in\CC,$$
where $\Gamma$ denotes the Gamma function.
For the construction of characteristic functions in sectors $S_\alpha$ for $\alpha\in(0,1]$ we will take $A=2-\alpha$ and $B=4-\alpha$ and we set
\begin{equation*}%\label{goodEalpha}
	\widetilde{E}_{\alpha}(z):=E_{2-\alpha, 4-\alpha}(-z){=\sum_{j=0}^\infty \frac{(-1)^jz^j}{\Gamma ((2-\alpha)(j+1)+2)},\quad z\in\CC}.
\end{equation*}
%For the construction of characteristic functions in sectors $S_\alpha$ for $\alpha\in[0,2)$ we will take $A=2-\alpha$ and $B=4-\alpha$ and so we set
%\begin{equation}\label{goodEalpha}
%	\widetilde{E}_{\alpha}(z):=E_{2-\alpha, 4-\alpha}(-z).
%\end{equation}

We recall the following statement.

\begin{theorem}\label{Ealphathm} (\cite[Thm. 5, Thm. 20]{salinas62})
%	Let $\widetilde{E}_{\alpha}(z)$ be the function from \eqref{goodEalpha}.
%\red{CASE $\alpha=0$.}
%	\begin{itemize}
%		\item[$(i)$] By \cite[Thm. 5]{salinas62} we get: Let $\alpha\in[0,2)$, then
%		$$\exists\;C\ge 1\;\forall\;z\in S_{\alpha}\;\forall\;n\in\NN_0:\;\;\;\left| \widetilde{E}_{\alpha}(z)-\sum_{j=0}^{n-1} \frac{(-z)^j}{\Gamma ((j+1) (2-\alpha)+2)}\right|\leq C \frac{e^n}{n^{(2-\alpha)n}}|z|^n.$$
%		Consequently, $\widetilde{E}_{\alpha}(z)\in\widetilde{\mathcal{A}}^u_{\{\overline{\G}^{\alpha-2}\}}(S_\alpha)$ and, moreover, $\widetilde{E}_{\alpha}(z)$ is a characteristic function in the class $\widetilde{\mathcal{A}}^u_{\{\overline{\G}^{\alpha-2}\}}(S_\alpha)$.
%		
%		\item[$(ii)$]
Let $\alpha\in(0,1]$, then
		\begin{equation}\label{Ealphagoodestimate}
			\forall\;z\in S_{\alpha}\;\forall\;n\in\NN_0:\;\;\;\left|\widetilde{E}_{\alpha}^{(n)}(z)\right|\leq 2 \frac{n! e^n}{n^{(2-\alpha)n}}.
		\end{equation}
		Consequently, $\widetilde{E}_{\alpha}\in \mathcal{A}_{\{\overline{\G}^{\alpha-1}\}}(S_\alpha)$.
		Moreover,
{$$
\widetilde{E}_{\alpha}^{(n)}(0)=\frac{(-1)^n n!}{\Gamma ((2-\alpha)(n+1)+2)},\quad n\in\NN_0,
$$
and so} $\widetilde{E}_{\alpha}$ is a characteristic function in the class $\mathcal{A}_{\{\overline{\G}^{\alpha-1}\}}(S_\alpha)$.
%	\end{itemize}
\end{theorem}

%Next we consider the \textbf{Rodr\'{i}guez-Salinas function} (see \cite{salinas62})
%\begin{equation}\label{RSfunction}
%	K_{RS}(z):=\frac{(1+z)\log(1+z)-z}{z^2},
%\end{equation}
%which is clearly holomorphic in $S_2$.
%
%Now we recall \cite[Thm. 3]{salinas62}
%
%\begin{theorem}\label{RStheorem}
%	There exists a constant $C$ (we can take $C=4$) such that
%	 $$\forall\;z\in\CC\backslash(-\infty,-1]\;\forall\;n\in\NN_0:\;\;\;\left| K_{RS}(z)-\sum_{j=0}^{n-1} \frac{(-z)^j}{(j+1)(j+2)}\right|\leq C|z|^n.$$
%	Consequently, for the constant sequence $\mathbbm{1}:=(1)_{j\in\NN_0}$ we have that  $K_{RS}\in\widetilde{\mathcal{A}}^u_{\{\mathbbm{1}\}}(S_2)$ and, moreover, $K_{RS}$ is a characteristic function in the class $\widetilde{\mathcal{A}}^u_{\{ \mathbbm{1}\}}(S_2)$.
%\end{theorem}
%
%
%Let $\alpha>2$ and take $\alpha'>\alpha$. For all $z\in S_\alpha$ we set
%\begin{equation}\label{fctfalphaalphaprime}
%	f_{\alpha,\alpha'}(z):= \int_{0}^{\infty(-\phi)} K_{RS}(zv^{\alpha'-2}) e^{-v} dv,
%\end{equation}
%where we choose $\phi  \in (-\frac{(\alpha-2)}{(\alpha'-2)}\frac{\pi}{2}, \frac{(\alpha-2)}{(\alpha'-2)}\frac{\pi}{2})$ with $|\arg(z)- (\alpha'-2)\phi|<\pi$.

Let $\alpha>1$ and take $\alpha'>\alpha$. For all $z\in S_\alpha$ we define
\begin{equation}\label{fctgalphaalphaprime}
	g_{\alpha,\alpha'} (z):= \int_{0}^{\infty(-\phi)} e^{-zv^{\alpha'-1}} e^{-v} dv,
\end{equation}
where we choose $\phi  \in (-\frac{(\alpha-1)}{(\alpha'-1)}\frac{\pi}{2}, \frac{(\alpha-1)}{(\alpha'-1)}\frac{\pi}{2})$ with $|\arg(z)- (\alpha'-1)\phi|<\pi/2$.

\begin{theorem}\label{fgalphaalphaprimethm} (\cite[Thm. 28]{salinas62})
	Let $\alpha>1$, $\alpha'>\alpha$ and
%$f_{\alpha,\alpha'}$ resp.
$g_{\alpha,\alpha'}$ be the function from
%\eqref{fctfalphaalphaprime} resp.
\eqref{fctgalphaalphaprime}.
%	\begin{itemize}
%		\item[$(i)$] By \cite[Thm. 4, Thm. 8]{salinas62} we get: Let $\alpha>2$ and take $\alpha'>\alpha$. Then
%		\begin{align*}
%			\exists\;C,A\ge 1\;\forall\;z\in S_{\alpha}\;\forall\;n\in\NN_0:\;\;\;&\left|f_{\alpha,\alpha'} (z) -\sum_{j=0}^{n-1} \frac{(-z)^j}{(j+1)(j+2)} \Gamma (  (\alpha'-2)j+1)\right|
%			\\&
%			\le CA^n \Gamma ( (\alpha'-2)n+1)|z|^n.
%		\end{align*}
%		Consequently, we have that $f_{\alpha,\alpha'}\in\widetilde{\mathcal{A}}^u_{\{\overline{\G}^{\alpha'-2}\}}(S_\alpha)$ and,
%		moreover, $f_{\alpha,\alpha'}$ is a characteristic function in the class $\widetilde{\mathcal{A}}^u_{\{ \overline{\G}^{\alpha'-2} \}}(S_\alpha)$.
		
%		\item[$(ii)$]
Then,
		\begin{equation}\label{galphagoodestimate}
			\exists\;C,A\ge 1\;\forall\;z\in S_{\alpha}\;\forall\;n\in\NN_0:\;\;\;\left| g_{\alpha,\alpha'}^{(n)}(z)\right|\leq C A^n \Gamma ((\alpha'-1)n+1).
		\end{equation}
		Consequently, $g_{\alpha,\alpha'}\in \mathcal{A}_{\{\overline{\G}^{\alpha'-1}\}}(S_\alpha)$. Moreover,
{$$
g_{\alpha,\alpha'}^{(n)}(0)=(-1)^n  \Gamma((\alpha'-1)n+1),\qquad n\in\NN_0,
$$
and so} $g_{\alpha,\alpha'}$ is a characteristic function of the class $\mathcal{A}_{\{\overline{\G}^{\alpha'-1}\}}(S_\alpha)$.
%\end{itemize}
\end{theorem}

\subsection{Characteristic transform}\label{ssect.CharactTransf}

Following again the work of Rodríguez Salinas~\cite{salinas62}, we present a functional transform that modifies the derivatives at 0 of a function in a ultraholomorphic class with a precise control, what allows for the construction of characteristic functions in more general classes than the Gevrey ones, considered previously.

\begin{definition}
	Let $\M$ be an (lc) sequence, $\LL\in\RR_{>0}^{\NN_0}$, $S$ a sector and %$f\in\widetilde{\mathcal{A}}^u_{\{\LL\}} (S)$ resp.
$f\in \mathcal{A}_{\{\LL\}} (S)$. Then we define the $\mathcal{T}_{\M}-$transform of $f$ by
	\begin{equation*}%\label{Ttrasform}
		\mathcal{T}_{\M}(f)(z):= \sum_{j=0}^{\infty}\frac{1}{2^j} \frac{M_j}{m_j^j} f(m_j z), \qquad z\in S.
	\end{equation*}
\end{definition}
{This expression} should be compared with the characteristic functions obtained in the ultradifferentiable setting in \cite[Thm. 1]{thilliez} and \cite[Lemma 2.9]{compositionpaper}. For every $j\in\NN_0$ let us set
\begin{equation*}%\label{Sexpression}
	R_j:=\sum^{\infty}_{n=0} \frac{1}{2^n} \frac{M_n}{m_n^n}m_n^j.
\end{equation*}
The following result provides estimates for this sequence in terms of the general sequence $\M$ we depart from.

\begin{lemma}\label{Scomplemma}
	Let $\M\in\RR_{>0}^{\NN_0}$, then
	$$\forall\;j\in\NN_0:\;\;\;R_j\ge\frac{1}{2^j} M_j.$$
	If $\M$ is (lc), then also
	$$\forall\;j\in\NN_0:\;\;\;R_j\le 2 M_j,$$
	and so $(R_j)_{j\in\NN_0}$ is equivalent to $\M$.
\end{lemma}

\demo{Proof}
For any $j\in\NN_0$ we choose $n=j$ in the sum and get $R_j\ge\frac{1}{2^j} \frac{M_j}{m_j^j}m_j^j=\frac{1}{2^j} M_j$.

For the converse we recall that since $\M$ is (lc) we have $m_0\le m_1\le\dots$ and so
$$\forall\;j,n\in\NN_0:\;\;\;(m_n)^{j-n}\le\frac{M_j}{M_n},$$
see \cite[Thm. 1]{thilliez} and the detailed proof in \cite[(3.1.2)]{diploma}. Thus
$$R_j=\sum^{\infty}_{n=0}\frac{1}{2^n}M_n m_n^{j-n}\le\sum^{\infty}_{n=0}\frac{1}{2^n}M_n\frac{M_j}{M_n}=2M_j$$
for all $j\in\NN_0$.
\qed\enddemo

\begin{theorem}\label{CharacteristicTransform}
	%Let $\M$ be a (lc) sequence, $\LL\in\RR_{>0}^{\NN_0}$ and for each $\beta>0$ take
 Let $\M$ be a (lc) sequence, $\LL\in\RR_{>0}^{\NN_0}$ and for a given sector $S$ take $f\in \mathcal{A}_{\{\LL\}} (S)$.%
 %For all $n\in\NN_0$ put  $$b_n:=\lim_{z\in S, z\rightarrow 0} f^{(n)}(z),$$ and
 Then,
%  we have:
%	\begin{enumerate}
%		\item $\mathcal{T}_{\M} (f)$ is holomorphic in $S$.
%		\item For every $j\in\NN_0$ we define
%		$$a_j:=(f^{(j)}(0)) R_j,$$
%		with $R_j$ being the numbers from \eqref{Sexpression}. Then
$\mathcal{T}_{\M} (f)\in \mathcal{A}_{\{\LL\M\}} (S)$ with
\begin{equation}\label{eq.DerivAt0T_Mf}
\mathcal{T}_{\M} (f)^{(j)}(0)=R_jf^{(j)}(0),\qquad j\in\NN_0.
\end{equation}
Moreover, for any $A>0$, $\mathcal{T}_{\M} : \mathcal{A} _{\LL, A} (S)\rightarrow\mathcal{A}_{\LL\M, A } (S)$ is a continuous linear operator.
%	\end{enumerate}
\end{theorem}

\demo{Proof}
%\begin{enumerate}
%	\item
By definition of $\mathcal{A}_{\{\LL\}} (S)$ we have that $f$ is bounded in $S$ by some constant $C>0$. Since $\M$ is log-convex, we have that $M_j\le m_j^j$ for all $j\in\NN_0$ and then
	$$\sum_{j=0}^{\infty}\frac{1}{2^j} \frac{M_j}{m_j^j} \left|f(m_j z) \right| \leq  C  \sum_{j=0}^{\infty}\frac{1}{2^j}=2C, \qquad z\in S.$$
Consequently, the series defining $\mathcal{T}_{\M} (f)$ normally converges in the whole of $S$, it provides a function holomorphic in $S$, and differentiation and limits can be interchanged with summation.
 %and so the proof follows as in Theorem \ref{CharacteristicTransformUniformAsymptotic}.
%	\item
For each $z\in S$ and every $j\in\NN_0$ we observe then that
	$$(\mathcal{T}_{\M} (f))^{(j)} (z)= \sum_{n=0}^{\infty}\frac{1}{2^n} \frac{M_n}{m_n^n} m_n^jf^{(j)}(m_n z),$$
and so
$$
\mathcal{T}_{\M} (f)^{(j)}(0)=\sum_{n=0}^{\infty}\frac{1}{2^n} \frac{M_n}{m_n^n} m_n^jf^{(j)}(0)=R_jf^{(j)}(0),\qquad j\in\NN_0,$$
as desired.

Suppose $f\in \mathcal{A}_{\LL,A} (S)$ for some $A>0$, then for all $j\in\NN_0$ we can estimate
\begin{align*}
|(\mathcal{T}_{\M} (f))^{(j)} (z)|&\le \sum_{n=0}^{\infty}\frac{1}{2^n} \frac{M_n}{m_n^n}m_n^j |f^{(j)}(m_nz)|
 \\&
 \le \|f\|_{\M,A} A^j L_j \sum_{n=0}^{\infty}\frac{1}{2^n}M_nm_n^{j-n}=\|f\|_{\M,A} A^j L_j R_j.
\end{align*}
By Lemma~\ref{Scomplemma} we know that $R_j\le 2M_j$, so $\mathcal{T}_{\M} (f)\in \mathcal{A}_{\LL\M, A } (S)$, and moreover
$$\|\mathcal{T}_{\M} (f)\|_{\LL\M, A}=\sup_{z\in S}\frac{|(\mathcal{T}_{\M} (f))^{(j)} (z)|}{A^j L_j M_j}\le 2\|f\|_{\M,A}.$$
It follows that $\mathcal{T}_{\M} : \mathcal{A} _{\LL, A} (S)\rightarrow\mathcal{A}_{\LL\M, A } (S)$ is a well-defined continuous linear operator for any $A>0$.
%\end{enumerate}
\qed\enddemo

\begin{theorem}\label{CharacteristicTransform1}
	%Let $\M$ be a (lc) sequence, $\LL\in\RR_{>0}^{\NN_0}$ and for each $\beta>0$ take
 Let $\M$ be a (lc) sequence, $\LL\in\RR_{>0}^{\NN_0}$ and for a given sector $S$ take $f\in \mathcal{A}_{\{\LL\}} (S)$. %
 %For all $n\in\NN_0$ we denote  $$b_n:=\lim_{z\in S, z\rightarrow 0} f^{(n)}(z).$$
	If $(|f^{(j)}(0)|)_{j\in\NN_0}$ is equivalent to $\LL$, then
$(|\mathcal{T}_{\M} (f)^{(j)}(0)|)_{j\in\NN_0}$
%	$$ (f^{(j)}(0)) R_j:=\mathcal{T}_{\M} (f)^{(j)}(0),$$
	is equivalent to $\LL\M$.
	Consequently, $\mathcal{T}_{\M} (f)$ is characteristic in the class $\mathcal{A}_{\{\LL\M\}} (S)$.
\end{theorem}

\demo{Proof}
The first assertion is clear from Lemma~\ref{Scomplemma} and~\eqref{eq.DerivAt0T_Mf}. The second one stems from Theorem~\ref{optthm2}.
%Since $(|f^{(j)}(0)|)_{j\in\NN_0}$ is equivalent to $\LL$ and by Lemma \ref{Scomplemma} we get that
%$$\exists\;C,A\le 1\;\forall\;n\in\NN_0:\;\;\;|(f^{(j)}(0)) R_n |\geq  C A^n L_n M_n,$$
%and the converse estimate follows analogously. Thus $(|(f^{(j)}(0))R_j|)_{j\in\NN_0}$ is equivalent to $\LL\M$.
\qed\enddemo

\subsection{Construction of characteristic functions}\label{ssect.ConstrCharactFunct}
Given a sequence $\M\in\RR_{>0}^{\NN_0}$ and $\alpha>0$ we construct now, under suitable assumptions, characteristic functions in $\mathcal{A}_{\{\M\}} (S_{\alpha})$.
% and in $\widetilde{\mathcal{A}}^u_{\{\M\}} (S_{\alpha})$.
For this we are using the basic functions from Subsection \ref{ssect.BasicFct} and the characteristic transform from Subsection \ref{ssect.CharactTransf}.\vspace{6pt}

\begin{theorem}\label{mainTthm2}
	Let $\M\in\RR_{>0}^{\NN_0}$ and $\alpha>0$.
%We assume that  $\overline{\G}^{1-\alpha'}\M:=(j^{(1-\alpha')j}M_j)_{j\in\NN_0}$
%	is equivalent to an (lc) sequence $\LL$ depending on $\alpha'$, where $\alpha'=\alpha$, if $\alpha\le 1$, or $\alpha'>\alpha$, if $\alpha>1$. Then, the following assertions hold true:
	\begin{enumerate}
		\item If $\alpha\le 1$, we assume that  $\overline{\G}^{1-\alpha}\M:=(j^{(1-\alpha)j}M_j)_{j\in\NN_0}$
	is equivalent to an (lc) sequence $\LL$. Then, $\mathcal{T}_{\LL} (\widetilde{E}_{\alpha}) $ is characteristic in the class $\mathcal{A}_{\{\M\}} (S_{\alpha})$.
		\item If $\alpha>1$, we assume that there exists $\a'>\a$ such that $\overline{\G}^{1-\alpha'}\M:=(j^{(1-\alpha')j}M_j)_{j\in\NN_0}$
	is equivalent to an (lc) sequence $\LL$. Then, $\mathcal{T}_{\LL} (g_{\alpha,\alpha'}) $ is characteristic in the class $\mathcal{A}_{\{\M\}} (S_{\alpha})$.
	\end{enumerate}
\end{theorem}

\demo{Proof}
This follows by Theorems~\ref{Ealphathm}, \ref{fgalphaalphaprimethm},  \ref{CharacteristicTransform} and \ref{CharacteristicTransform1}, and  from the fact that $\overline{\G}^{\alpha-1}\LL$ in case 1, resp. $\overline{\G}^{\alpha'-1}\LL$ in case 2, is equivalent to $\M$.
\qed\enddemo

\begin{remark}\label{rem.gammaindexcharactfunct}
In order to guarantee that the hypotheses in the previous theorem are satisfied, one can compute the index $\gamma(\M)$ and check whether it is greater than $\a-1$. If this is the case, the very definition of this index implies that for {any $\beta$ such that $\gamma(\M)>\beta>\a-1$ the property $\left(P_{\beta}\right)$} (see Subsection~\ref{ssect.IndexgammaM}) is satisfied, and so there exists a suitable  (lc) sequence $\LL$ in the desired conditions.
\end{remark}

%\red{Theorem: With the previous notations we get:
%	\begin{enumerate}
%		\item If $\alpha\le 1$, then $\mathcal{T}_{\overline{\G}^{1-\alpha}\M} (E_\alpha) $ is characteristic in the class $\mathcal{A}_{\{\M\}} (S_{\alpha})$.
%		\item If $\alpha>1$, then $\mathcal{T}_{\overline{\G}^{1-\alpha'}\M} (g_{\alpha,\alpha'}) $ is characteristic in the class $\mathcal{A}_{\{\M\}} (S_{\alpha})$.
%\end{enumerate}}

\section{Stability properties for ultraholomorphic classes defined by weight matrices}\label{sect.StabilityPropert}

The aim of this section is to generalize and extend the stability result of Ider and Siddiqi~\cite[Thm. 1]{IderSiddiqi}, valid for Carleman-Roumieu ultraholomorphic classes in sectors not wider than a half-plane. We give the proof in the general weight matrix setting, we get rid of the restriction on the opening of the sector (thanks to the construction of characteristic functions in arbitrary sectors), and we extend the list of stability properties.

%\subsection{Basic notation}

%\subsection{Main result}
Our main result is concerned with several stability properties which will be defined next.

\begin{definition} Let $\M\in\RR_{>0}^{\NN_0}$ be a sequence and $U\subseteq\CC$ be an open set. Given a compact set $K\subset U$, we define
	$$\mathcal{H}_{\M,h}(K):=\{f\in \mathcal{H}(U): \|f\|_{\M,K,h}:=\sup_{z\in K, j\in\NN_0}\frac{|f^{(j)}(z)|}{h^j M_j}<+\infty\}.$$
We put
	\begin{equation*}
		 \mathcal{H}_{\{\M\}}(K):=\bigcup_{h>0}\mathcal{H}_{\M,h}(K).
	\end{equation*}
Moreover, given a weight matrix $\mathcal{M}=\{\M^{(p)}: p>0\}$, we may introduce the class $\mathcal{H}_{\{\mathcal{M}\}}(U)$ as
	\begin{equation*}%\label{generalroumieuultraholom}
		\mathcal{H}_{\{\mathcal{M}\}}(U):=\bigcap_{K\subset U}\bigcup_{p>0}\mathcal{H}_{\{\M^{(p)}\}}(K).
	\end{equation*}	
\end{definition}

\begin{definition} Let $\mathcal{M}=\{\M^{(p)}: p>0\}$ be a weight matrix and $\alpha>0$. The class $\mathcal{A}_{\{\mathcal{M}\}}(S_{\alpha})$ is said to be:
		\begin{itemize}
			\item[$(i)$] holomorphically closed, if for all $f\in\mathcal{A}_{\{\mathcal{M}\}}(S_{\alpha})$ and $g\in\mathcal{H}(U)$, where $U\subseteq\CC$ is an open set containing the closure of the range of $f$, we have $g\circ f\in\mathcal{A}_{\{\mathcal{M}\}}(S_{\alpha})$.
			\item[$(ii)$] inverse-closed, if for all $f\in\mathcal{A}_{\{\mathcal{M}\}}(S_{\alpha})$ such that  $\inf_{z\in S_\alpha}|f(z)|>0$, we have $1/f\in\mathcal{A}_{\{\mathcal{M}\}}(S_{\alpha})$.
			\item[$(iii)$] closed under composition, if for all $f\in\mathcal{A}_{\{\mathcal{M}\}}(S_{\alpha})$ and for all $g\in\mathcal{H}_{\{\mathcal{M}\}}(U)$, where $U\subseteq\CC$ is an open set containing the closure of the range of $f$, we have $g\circ f\in\mathcal{A}_{\{\mathcal{M}\}}(S_{\alpha})$.
		\end{itemize}
\end{definition}

\begin{remark}
We wish to highlight that it is important to state these definitions in a clear way. We cannot relax the condition $\inf_{z\in S_\alpha}|f(z)|>0$ in the definition of inverse-closed by considering, for example, the weaker requirement{:
$$f(z)\neq0\qquad \text{for all $ z\in S_\alpha $.}$$} While this is enough when working with ultradifferentiable classes on compact intervals, as done in~\cite{Malliavin}, our situation is different as ${S_\alpha}$ is not compact. This is easily seen by considering the function $z\mapsto \exp(-1/z)$, which belongs to the class $ \mathcal{A}_{\{\G^2\}}(S_{\a}) $ for every $\a\in(0,1)$ {(as a consequence of Cauchy's integral formula for the derivatives)} and never vanishes in $S_\a$. However, observe that its multiplicative inverse $z\mapsto \exp(1/z)$ is not bounded, and hence it does not belong to any of the ultraholomorphic classes under consideration.

In the same vein, the open set $U$ in (i) and (iii) has to contain the closure of the range of $f$, and not just the range. This is clearly seen in the forthcoming arguments involving the function $z\mapsto 1/z$, whose derivatives admit global analytic bounds in closed subsets of $\CC\setminus\{0\}$, but not in the whole of it.
\end{remark}

Our first statement will consider classes in sectors $S_\a$ contained in a half-plane and defined by a weight matrix $\mathcal{M}$. In this case, the matrix can be changed, without altering the class, into a new matrix $\mathcal{M}^\a$ which we define now.

{
\begin{definition}
Let $\mathcal{M}=\{\M^{(p)}: p>0\}$ be a weight matrix (not necessarily satisfying \hyperlink{Msc}{$(\mathcal{M}_{\on{sc}})$}). Given $\alpha>0$ we assume that $\lim_{j\rightarrow+\infty}(j^{(1-\alpha)j}M^{(p)}_j)^{1/j}=\infty$ for all $p>0$. The matrix
\begin{equation*}%\label{matrixalpha}
	\mathcal{M}^{\alpha}:=\{\M^{(p,\alpha)}: p>0\}
\end{equation*}
is defined as
\begin{equation}\label{ramiequ}
\M^{(p,\alpha)}=\overline{\G}^{\alpha-1}\left(\overline{\G}^{1-\alpha}\M^{(p)}\right)^{\on{lc}},
\hspace{15pt}
M^{(p,\alpha)}_j=j^{(\alpha-1)j}\left(\overline{\G}^{1-\alpha}\M^{(p)}\right)^{\on{lc}}_j,
%\left(j^{(1-\alpha)j}M^{(p)}_j\right)^{\on{lc}}=j^{(1-\alpha)j},
\;\;\;\;j\in\NN_0.
\end{equation}
\end{definition}
}

So, every sequence in the original matrix is termwise multiplied by the Gevrey-like sequence $\overline{\G}^{1-\alpha}$ (recall that $\overline{\G}^{1-\alpha}\hyperlink{approx}{\approx}\G^{1-\alpha}$), this sequence is changed into its log-convex regularization, and finally one termwise divides by $\overline{\G}^{1-\alpha}$ again.
It is clear that $M^{(p,\alpha)}_0=M^{(p)}_0=1$ (recall the convention $0^0:=1$) for all $\alpha>0$ and $p>0$, and that the map $p\mapsto M^{(p,\alpha)}_j$ is non-decreasing for any $j\in\NN_0$ fixed.
So, $\M^{(p,\alpha)}\le \M^{(p',\alpha)}$ for all $0<p< p'$, i.e., $\mathcal{M}^{\alpha}$ is a weight matrix according to the definition given in Subsection \ref{ssect.WeightMatrices}. However, in general $\mathcal{M}^{\alpha}$ is not log-convex.
%{\color{red}Finally let us note that the technical switch $\M^{(p)}\mapsto(\M^{(p)})^{\on{lc}}$ as it has been done in \cite{IderSiddiqi} is superfluous since by our general assumption $\mathcal{M}$ is \hyperlink{Msc}{$(\mathcal{M}_{\on{sc}})$}.}
%\begin{remark} Probably we only need this for all $p>0$: $$\liminf_{j\rightarrow+\infty}\frac{(M^{(p)}_j)^{1/j}}{n^{(1-\alpha)j}}>0.$$
%	Please check Thm 20 and the footnote (19) after the Thm 24 in Salinas (1962). Maybe if this $\liminf$ is finite, the corresponding class is equivalent to the class $\mathcal{A}_{\{\mathbbm{1}\}}(S_{\alpha})$. In any case, is important to comment that Sidiqqi should write this hypothesis in his theorem.
%\end{remark}

{
\begin{remark}\label{liminf-remark} Note that if there exist some $p>0$ such that $\lim_{j\rightarrow+\infty}(j^{(1-\alpha)j}M^{(p)}_j)^{1/j}=\infty$, then the same is valid for all $p'>p$, thanks to the fact that the $\M^{(p)}\leq \M^{(p')}$. In this situation, since we also have $\mathcal{A}_{\{\M^{(p)}\}}(S_{\alpha})\subseteq \mathcal{A}_{\{\M^{(p')}\}}(S_{\alpha})$ and the class associated to the weight matrix $\mathcal{M}$ is the increasing union of such classes, in order to study stability properties in it we can restrict our attention to the case described in the previous definition.

In case
$\lim_{j\rightarrow+\infty}(j^{(1-\alpha)j}M^{(p)}_j)^{1/j}$ is not infinity for any $p>0$, then there are some possibilities:
\begin{enumerate}
	\item[(i)] If $\alpha>1$ and $\liminf_{j\rightarrow+\infty}(j^{(1-\alpha)j}M^{(p)}_j)^{1/j}<\infty$ for all $p>0$, the class $\mathcal{A}_{\{\M^{(p)}\}}(S_{\alpha})$ only contains constant functions, see \cite[Thm. 21, and p. 8]{salinas62}, and the same holds for the class $\mathcal{A}_{\{\mathcal{M}\}}(S_{\alpha})$. So, the stability results turn out to be trivial.
	
	\item[(ii)] If $0<\alpha\le 1$ and $\liminf_{j\rightarrow+\infty}(j^{(1-\alpha)j}M^{(p)}_j)^{1/j}=0$ for all $p>0$, the class $\mathcal{A}_{\{\M^{(p)}\}}(S_{\alpha})$ only contains constant functions, see \cite[Thm. 20]{salinas62}, and again we are done.
	
	\item[(iii)] If $0<\alpha\le 1$ and $\liminf_{j\rightarrow+\infty}(j^{(1-\alpha)j}M^{(p)}_j)^{1/j}\in(0,\infty)$ for all $p>0$ (or from some $p_0>0$ on), taking into account \cite[Cor. 8]{salinas62} we have that the class $\mathcal{A}_{\{\M^{(p)}\}}(S_{\alpha})$ coincides with $\mathcal{A}_{\{\overline{\G}^{\alpha-1}\}}(S_{\alpha})$ for all $p>0$ (or for $p\ge p_0$), and so $\mathcal{A}_{\{\mathcal{M}\}}(S_{\alpha})= \mathcal{A}_{\{\mathcal{\overline{G}}^{\alpha-1}\}}(S_{\alpha})$, where $\mathcal{\overline{G}}^{\alpha-1}$ is the matrix with all the rows equal to the sequence $\overline{\G}^{\alpha-1}$. We will study the stability properties for this class in Section~\ref{sect.Examples}.
\end{enumerate}
%
%
%
%Please check Thm 20 and the footnote (19) after the Thm 24 in Salinas (1962). Maybe if this $\liminf$ is finite, the corresponding class is equivalent to the class $\mathcal{A}_{\{\mathbbm{1}\}}(S_{\alpha})$. In any case, is important to comment that Sidiqqi should write this hypothesis in his theorem.
\end{remark}
}

In order to prove the aforementioned equality of the classes associated with $\mathcal{M}$ and $\mathcal{M}^{\a}$, it is convenient to recall the following result, {which provides Gorny-Cartan like inequalities for holomorphic functions in sectors.}

\begin{theorem}\label{Gornycartan} (\cite[Thm. 23]{salinas62}) Let $0<\alpha\leq1$ and $f\in \mathcal{H}(S_\alpha)$. If
	$C_n(f)=\sup_{z\in S_\alpha}|f^{(n)}(z)|$, $n\in\NN_0$,
	then the sequence $B_n=n^{(1-\alpha)n}C_n(f)$ verifies
	\begin{equation*}%\label{Gornycartaninequality}
		B_n\leq A q^{(1-\alpha)n}B_{n_1}^{\frac{n_2-n}{n_2-n_1}}B_{n_2}^{\frac{n-n_1}{n_2-n_1}},\qquad n_1< n < n_2,
	\end{equation*}
	where $A=4$ and $q=1$ if $\alpha=1$, or $A=8\pi$ and $q=2e(2-\alpha)/(1-\alpha)$ for the remaining cases.
\end{theorem}

\begin{theorem}\label{classequality} Let $\mathcal{M}=\{\M^{(p)}: p>0\}$ be a weight matrix and $0<\alpha\le 1$ be given such that $\lim_{j\rightarrow+\infty}(j^{(1-\alpha)j}M^{(p)}_j)^{1/j}=\infty$ for all $p>0$. Let $\mathcal{M}^\alpha=\{\M^{(p,\alpha)}: p>0\}$ be the matrix given in~\eqref{ramiequ}. Then, we have that
$$
\mathcal{A}_{\{\mathcal{M}\}}(S_{\alpha})= \mathcal{A}_{\{\mathcal{M}^\alpha\}}(S_{\alpha}).
$$	
\end{theorem}
\demo{Proof}
Given $f\in\mathcal{A}_{\{\mathcal{M}^\alpha\}}(S_{\alpha})$, there exists some $p>0$ such that $f\in \mathcal{A}_{\{\M^{(p,\alpha)}\}}(S_{\alpha})$. Since $\overline{\G}^{1-\alpha}\M^{(p,\alpha)}$ is the log convex minorant of $\overline{\G}^{1-\alpha}\M^{(p)}$, we have that $\overline{\G}^{1-\alpha}\M^{(p,\alpha)}\leq\overline{\G}^{1-\alpha}\M^{(p)}$, and therefore $\M^{(p,\alpha)}\leq\M^{(p)}$. We conclude that $f\in \mathcal{A}_{\{\mathcal{M}\}}(S_{\alpha})$.\\
For the converse inclusion, let us consider $f\in\mathcal{A}_{\{\mathcal{M}\}}(S_{\alpha})$. There exist some $C,D\in\RR_{>0}$ and $p>0$ such that $C_n(f)=\sup_{z\in S_\alpha}|f^{(n)}(z)|\leq C D^n M_n^{(p)}$, for all $n\in\NN_0$.

Let us fix $n\in\NN_{0}$ and distinguish two cases:
\begin{enumerate}
	\item [i)] If $M_n^{(p,\alpha)}=M_n^{(p)}$ then $\sup_{z\in S_\alpha}|f^{(n)}(z)|\leq C D^n M_n^{(p,\alpha)}$.
	\item[ii)] If not, by the construction of
%thanks to the fact that $\overline{\G}^{1-\alpha}\M^{(p,\alpha)}$ is
the log convex minorant,
% of the sequence $\overline{\G}^{1-\alpha}\M^{(p)}$,
there exist so-called principal indices $n_1,n_2\in\NN_0$, with $n_1<n<n_2$, such that $M_{n_i}^{(p,\alpha)}=M_{n_i}^{(p)}$ for $i=1,2$ (see~\cite[Chapitre~I]{mandelbrojtbook} and, for a detailed discussion of the regularization process and its intricacies,~\cite{regularization}). So, we have
	\begin{align*}
		 \ln(n^{(1-\alpha)n}M_n^{(p,\alpha)})&=\frac{n_2-n}{n_2-n_1}\ln(n_1^{(1-\alpha)n_1}M_{n_1}^{(p,\alpha)}) +\frac{n-n_1}{n_2-n_1}\ln(n_2^{(1-\alpha)n_2}M_{n_2}^{(p,\alpha)})\\
		&\geq \frac{n_2-n}{n_2-n_1}\ln(\frac{1}{CD^{n_1}}n_1^{(1-\alpha)n_1}C_{n_1}(f))\\ &+\frac{n-n_1}{n_2-n_1}\ln(\frac{1}{CD^{n_2}}n_2^{(1-\alpha)n_2}C_{n_2}(f)).
		 %&=\frac{n_2-n}{n_2-n_1}\ln(\frac{1}{CD^{n_1}}B_{n_1})+\frac{n-n_1}{n_2-n_1}\ln(\frac{1}{CD^{n_2}}B_{n_2}).
	\end{align*}
	Therefore, with the notation of the previous theorem, we deduce from above:
	\begin{equation*}
		 B_{n_1}^{\frac{n_2-n}{n_2-n_1}}B_{n_2}^{\frac{n-n_1}{n_2-n_1}}\leq (CD^{n_1})^{\frac{n_2-n}{n_2-n_1}}(CD^{n_2})^{\frac{n-n_1}{n_2-n_1}}n^{(1-\alpha)n}M_n^{(p,\alpha)} =CD^nn^{(1-\alpha)n}M_n^{(p,\alpha)}.
	\end{equation*}
	Now, from the previous estimate and by applying Theorem~\ref{Gornycartan}, there exist some $A,q>0$ such that
	\begin{equation*}
		C_n(f)\leq n^{(\alpha-1)n}A q^{(1-\alpha)n}B_{n_1}^{\frac{n_2-n}{n_2-n_1}}B_{n_2}^{\frac{n-n_1}{n_2-n_1}}\leq AC (q^{(1-\alpha)}D)^nM_n^{(p,\alpha)}.
	\end{equation*}
\end{enumerate}
We conclude that $f\in \mathcal{A}_{\{\mathcal{M}^\alpha\}}(S_{\alpha})$.
\qed\enddemo

We are ready to state our first main result.

\begin{theorem}\label{theorem1siddiqi}
	Let $\mathcal{M}=\{\M^{(p)}: p>0\}$ be a weight matrix (not necessarily \hyperlink{Msc}{$(\mathcal{M}_{\on{sc}})$}) and $0<\alpha\le 1$ be given such that $\lim_{j\rightarrow+\infty}(j^{(1-\alpha)j}M^{(p)}_j)^{1/j}=\infty$ for all $p>0$. Let $\mathcal{M}^{\alpha}=\{\M^{(p,\alpha)}: p>0\}$ be the matrix according to \eqref{ramiequ}. Then the following assertions are equivalent:
	\begin{itemize}
		\item[$(a)$] The matrix $\mathcal{M}^{\alpha}$ satisfies the property \hyperlink{R-rai}{$(\mathcal{M}_{\{\on{rai}\}})$}.
		
		\item[$(b)$] The class $\mathcal{A}_{\{\mathcal{M}\}}(S_{\alpha})$ is holomorphically closed.

		\item[$(c)$] The class $\mathcal{A}_{\{\mathcal{M}
			\}}(S_{\alpha})$ is inverse-closed.
		
	\end{itemize}
	If $\mathcal{M}$ has in addition \hyperlink{R-Comega}{$(\mathcal{M}_{\{\text{C}^{\omega}\}})$} and $\mathcal{M}^{\alpha}$ has \hyperlink{R-dc}{$(\mathcal{M}_{\{\on{dc}\}})$}, then the list of equivalences can be extended by

	\begin{itemize}
		\item[$(d)$] The class $\mathcal{A}_{\{\mathcal{M}\}}(S_{\alpha})$ is closed under composition.

		\item[$(e)$] The matrix $\mathcal{M}^{\alpha}$ satisfies the property \hyperlink{R-FdB}{$(\mathcal{M}_{\{\on{FdB}\}})$}.
	\end{itemize}
\end{theorem}

\demo{Proof}
$(a)\Rightarrow(b)$ First recall that by the so-called {\itshape Fa\`{a}-di-Bruno formula} for the composition we get
$$(g\circ f)^{(n)}(z)=\sum_{\sum_{i=1}^nk_i=k, \sum_{i=1}^nik_i=n}\frac{n!}{k_1!\cdots k_n!}g^{(k)}(f(z))\prod_{i=1}^n\left(\frac{f^{(i)}(z)}{i!}\right)^{k_i},\;\;\;z\in S_{\alpha},\;n\in\NN_0.$$
%Let us now see that \hyperlink{R-rai}{$(\mathcal{M}_{\{\on{rai}\}})$} for $\mathcal{M}^{\alpha}$ yields
%\begin{equation}\label{theorem1siddiqiequ0}
%	\forall\;p>0\;\exists\;H\ge 1\;\exists\;p'(\ge p)\;\forall\;k\in\NN\;\forall\;j_k\in\NN_0:\;\;\;\check{M}^{(p,\alpha)}_{j_1}\cdots \check{M}^{(p,\alpha)}_{j_k}\le H^{j_1+\dots+j_k} \check{M}^{(p',\alpha)}_{j_1+\dots+j_k}.
%\end{equation}
%If $j_k\ge 1$ for all $k$ we estimate by
%$$\check{M}^{(p,\alpha)}_{j_1}\cdots \check{M}^{(p,\alpha)}_{j_k}\le H^{j_1}(\check{M}^{(p',\alpha)}_{j_1+\dots+j_k})^{\frac{j_1}{j_1+\dots +j_k}}\cdots H^{j_k}(\check{M}^{(p',\alpha)}_{j_1+\dots+j_k})^{\frac{j_k}{j_1+\dots +j_k}}=H^{j_1+\dots+j_k}\check{M}^{(p',\alpha)}_{j_1+\dots+j_k},$$
%and the remaining cases follow by $\check{M}^{(p,\alpha)}_0=M^{(p,\alpha)}_0=1$. Note that the indices $p$ and $p'$ are related by property \hyperlink{R-rai}{$(\mathcal{M}_{\{\on{rai}\}})$}.
%
%
%$$\check{M}^{(p,\alpha)}_{j_1}\cdots \check{M}^{(p,\alpha)}_{j_k}\le H^{j_1}(\check{M}^{(p',\alpha)}_{j_1+\dots+j_k})^{j_1/(j_1+\dots j_k)}\cdots H^{j_k}(\check{M}^{(p',\alpha)}_{j_1+\dots+j_k})^{j_k/(j_1+\dots+j_k)}=H^{j_1+\dots+j_k}\check{M}^{(p',\alpha)}_{j_1+\dots+j_k},$$
%
Let now $f\in\mathcal{A}_{\{\mathcal{M}\}}(S_{\alpha})$ be given. By Theorem \ref{classequality} we know that the classes $\mathcal{A}_{\{\mathcal{M}^{\alpha}\}}(S_{\alpha})$ and $ \mathcal{A}_{\{\mathcal{M}\}}(S_{\alpha})$ are equal, therefore $f\in\mathcal{A}_{\{\mathcal{M}^{\alpha}\}}(S_{\alpha})$. In particular, $f$ is bounded and thus any function $g$ which is analytic in a domain containing the (compact) closure of the range of $f$ satisfies
\begin{equation}\label{analyticbounded}
	\exists\;C_1,h_1\ge 1\;\forall\;k\in\NN_0\;\forall\;z\in S_{\alpha}:\;\;\;|g^{(k)}(f(z))|\le C_1h_1^kk!.
\end{equation}
By applying this and the fact that $f\in\mathcal{A}_{\{\mathcal{M^\alpha}\}}$, we estimate as follows for all $n\in\NN_0$ and $z\in S_{\alpha}$:
\begin{align*}
	|(g\circ f)^{(n)}(z)|&\le\sum_{\sum_{i=1}^nk_i=k, \sum_{i=1}^nik_i=n}\frac{n!}{k_1!\cdots k_n!}|g^{(k)}(f(z))|\prod_{i=1}^n\left|\frac{f^{(i)}(z)}{i!}\right|^{k_i}
	\\&
	 \le\sum_{\sum_{i=1}^nk_i=k, \sum_{i=1}^nik_i=n}\frac{n!}{k_1!\cdots k_n!}C_1h_1^kk!\prod_{i=1}^n\left(C_2h_2^i \check{M}^{(p,\alpha)}_i\right)^{k_i}
	\\&
	\le C_1\sum_{\sum_{i=1}^nk_i=k, \sum_{i=1}^nik_i=n}\frac{n!}{k_1!\cdots k_n!}h_1^kk!C_2^{k_1+\dots+k_n}h_2^{k_1+\dots+ nk_n}\prod_{i=1}^n(\check{M}^{(p,\alpha)}_i)^{k_i}
	\\&
	 \underbrace{\le}_{\eqref{theorem1siddiqiequ0}}C_1(C_2h_1h_2)^n\sum_{\sum_{i=1}^nk_i=k, \sum_{i=1}^nik_i=n}\frac{n!}{k_1!\cdots k_n!}k!\prod_{i=1}^n H_1^{ik_i}\check{M}^{(p',\alpha)}_{ik_i}
	\\&
	 \underbrace{\le}_{\eqref{theorem1siddiqiequ0}}C_1(H_1C_2h_1h_2)^n\sum_{\sum_{i=1}^nk_i=k, \sum_{i=1}^nik_i=n}\frac{n!}{k_1!\cdots k_n!}k!H_2^{k_1+\dots+nk_n}\check{M}^{(p'',\alpha)}_{k_1+\dots+nk_n}
	\\&
	 =C_1(H_1H_2C_2h_1h_2)^nM^{(p'',\alpha)}_n\sum_{\sum_{i=1}^nk_i=k, \sum_{i=1}^nik_i=n}\frac{k!}{k_1!\cdots k_n!}
	\\&
	\le C_1C_3(H_1H_2C_2C_4h_1h_2)^nM^{(p'',\alpha)}_n.
\end{align*}
For the estimates also note that $k\le n$ and w.l.o.g. $C_2,h_1,h_2,H_1\ge 1$. Moreover, we have that
$$
\sum_{\sum_{i=1}^nk_i=k, \sum_{i=1}^nik_i=n}\frac{k!}{k_1!\cdots k_n!}=2^{n-1},
$$
see \cite[Lemma 1.4.1]{Krantz02} or \cite[Prop. 2.1]{FernandezGalbis06}. Finally, by taking into account that the classes $\mathcal{A}_{\{\mathcal{M}^{\alpha}\}}(S_{\alpha})$ and $ \mathcal{A}_{\{\mathcal{M}\}}(S_{\alpha})$ are equal, then $g\circ f\in\mathcal{A}_{\{\mathcal{M}\}}(S_{\alpha})$ is verified.\vspace{6pt}
	
%Note that citation [7] in \cite{Siddiqi} is wrong.

$(b)\Rightarrow(c)$ This is obvious by taking $g: z\mapsto\frac{1}{z}$ since $g\in\mathcal{H}(\CC\backslash\{0\})$ and $\CC\backslash\{0\}$ contains the (compact) closure of the image of any element $f\in\mathcal{A}_{\{\mathcal{M}\}}(S_{\alpha})$ such that  $\inf_{z\in S_\alpha}|f(z)|>0$.\vspace{6pt}

$(c)\Rightarrow(a)$ We follow the ideas from \cite[Thm. 1]{IderSiddiqi} and apply the constructions from the previous section. First, recall that $\LL^{(p)}:=\overline{\G}^{1-\alpha}\M^{(p,\alpha)}{=(\overline{\G}^{1-\alpha}\M^{(p,\alpha)})^{\text{lc}}}$ is log-convex for any $p>0$, see \eqref{ramiequ}. Let $p>0$ be arbitrary but from now on fixed. According to Theorem \ref{mainTthm2} we put
\begin{equation*}%\label{theorem1siddiqiequ}
	f_{p}(z):=\mathcal{T}_{\LL^{(p)}} (\widetilde{E}_{\alpha})(z).
	 %f_{p}(z):=\mathcal{T}_{\overline{\G}^{1-\alpha}\M^{(p,\alpha)}} (\widetilde{E}_{\alpha})(z),\;\alpha\in[0,1].
\end{equation*}
By using \eqref{Ealphagoodestimate} and Lemma \ref{Scomplemma} we estimate as follows:

\begin{align*}
	|f^{(n)}_{p}(z)|&\le\sum_{k= 0}^\infty\frac{1}{2^k}L^{(p)}_k\frac{(\ell^{(p)}_{k})^n}{(\ell^{(p)}_{k})^k}|\widetilde{E}_{\alpha}^{(n)}(\ell^{(p)}_{k}z)|\le 4L^{(p)}_n\frac{n!e^n}{n^{(2-\alpha)n}}
	\\&
	=4M^{(p,\alpha)}_n\frac{n!e^n}{n^n}\le 4e^nM^{(p,\alpha)}_n,
\end{align*}
%\begin{align*}
%	|f^{(n)}_{p}(z)|&\le\sum_{k\ge 0}\frac{1}{2^k}{k^{(1-\alpha)k}M^{(p,\alpha)}_k}\frac{(\mu^{(p)}_{{k+1}})^n}{(\mu^{(p)}_{{k+1}})^k}|\widetilde{E}_{\alpha}^{(n)}(\mu^{(p)}_{{k+1}}z)|\le 4{n^{(1-\alpha)n}M^{(p,\alpha)}_n}\frac{n!e^n}{n^{({2-\alpha})n}}
%	\\&
%	=4M^{(p,\alpha)}_n\frac{n!e^n}{n^n}\le 4e^nM^{(p,\alpha)}_n,
%\end{align*}
for all $n\in\NN_0$ and $z\in S_{\alpha}$.
%{\color{red}(because ${L^{(p)}_k\frac{(\ell^{(p)}_{k})^n}{(\ell^{(p)}_{k})^k}}\le {n^{(1-\alpha)n}M^{(p,\alpha)}_n}$ for all $(n,k)\in\NN_0^2$ by log-convexity for ${\LL^{(p)}}$)}.%${\overline{\G}^{1-\alpha}\M^{(p,\alpha)}}$.
This estimate shows that $f_{p}\in\mathcal{A}_{\{\mathcal{M}^{\alpha}\}}(S_{\alpha})$ and, in particular when being applied to $n=0$, it yields $\sup_{z\in S_{\alpha}}|f_{p}(z)|\leq 4<+\infty$.

Set $R^{(p)}_n:=\sum_{k= 0}^\infty\frac{1}{2^k}L^{(p)}_k(\ell^{(p)}_{k})^{n-k}$ and so we get
\begin{equation}\label{theorem1siddiqiequ3}
	 \forall\;n\in\NN_0:\;\;\;f^{(n)}_{p}(0)=R^{(p)}_n\frac{n!(-1)^n}{{\Gamma((2-\alpha)(n+1)+2)}},
\end{equation}
and from Lemma \ref{Scomplemma}
\begin{equation}\label{theorem1siddiqiequ4}
	 \forall\;n\in\NN_0:\;\;\;R^{(p)}_n\geq\frac{n^{(1-\alpha)n}M^{(p,\alpha)}_n}{2^n}.
\end{equation}
Take $\lambda>4$ (note that in \cite[p. 349, line 5]{Siddiqi} there is a mistake, one should write $\lambda>C_0(f)M_0^{\alpha}$).
%[Comment:{In our case $\sup_{z\in S_{\alpha}}|f_{p}(z)|<C(f_{p})M_0^{\alpha}=C(f_{p})$, so we can write $\lambda>C(f_{p})$. In addition, we can write $\lambda>4$ if we want a uniform $\lambda$ value.}{ \color{red}!!Compare with $M_0^{\alpha}$ in the paper; it seems to be not good enough choice!!} {I think that there is a mistake in the paper, they should have written $C(f)M_0^{\alpha}$ instead of $M_0^{\alpha}$}]
Then, if we put $\widetilde{f}_{p}:=\lambda-f_{p}$, we have that $\widetilde{f}_p\in\mathcal{A}_{\{\mathcal{M}^{\alpha}\}}(S_{\alpha})$. Moreover, since $\inf_{z\in S_\alpha}|\widetilde{f}_{p}(z)|>0$ and $\mathcal{A}_{\{\mathcal{M}^{\alpha}\}}(S_{\alpha})(=\mathcal{A}_{\{\mathcal{M}\}}(S_{\alpha}))$ is assumed to be inverse-closed, we get that $z\mapsto\frac{1}{\widetilde{f}_{p}(z)}=\frac{1}{\lambda-f_{p}(z)}\in\mathcal{A}_{\{\mathcal{M}^{\alpha}\}}(S_{\alpha})$. We write $g:z\mapsto\frac{1}{\lambda-z}$,
%\red{the dependence on $p$ is justified because $\lambda$ is clearly depending on this chosen index {(this comment isn't necessary if $\lambda>4$.)}}
then by applying again the Fa\`a-di-Bruno-formula to the composition $g\circ f_{p}\in\mathcal{A}_{\{\mathcal{M}^{\alpha}\}}(S_{\alpha})$ {and thanks to the fact that $g^{(k)}(z)=\frac{k!}{(\lambda-z)^{k+1}}$ for all $k\in\NN_0$,} yields: For some $C,h>0$ and some index $p'>0$ (large) we get for all $n\in\NN_0$ that
\begin{align*}
	|(g\circ f_{p})^{(n)}(0)|&=\left|\sum_{\sum_{i=1}^nk_i=k, \sum_{i=1}^nik_i=n}\frac{n!}{k_1!\cdots k_n!}\frac{k!}{(\lambda-f_{p}(0))^{k+1}}\prod_{i=1}^n\left(\frac{f_{p}^{(i)}(0)}{i!}\right)^{k_i}\right|\\
&\le Ch^nM^{(p',\alpha)}_n.
\end{align*}
%\red{Note that $g^{(k)}(z)=\frac{k!}{(\lambda-z)^{k+1}}$ for all $k\in\NN_0$.}
By \eqref{theorem1siddiqiequ3} we see
$$
\left(\frac{f_{p}^{(i)}(0)}{i!}\right)^{k_i}= \left(\frac{(-1)^iR^{(p)}_i}{{\Gamma((2-\alpha)(i+1)+2)}}\right)^{k_i}, \quad 1\le i\le n,
$$
and by taking into account that $ \prod_{i=1}^n (-1)^{ik_i}=(-1)^n $, we deduce that for every $n\in\NN_0$,
\begin{align*}
\sum_{\sum_{i=1}^nk_i=k, \sum_{i=1}^nik_i=n}\frac{n!}{k_1!\cdots k_n!}\frac{k!}{(\lambda-f_{p}(0))^{k+1}} \prod_{i=1}^n\left(\frac{R^{(p)}_i}{{\Gamma((2-\alpha)(i+1)+2)}}\right)^{k_i}
\le Ch^nM^{(p',\alpha)}_n.
\end{align*}
Each summand in this sum is strictly positive and we focus now on the one given by the choices $k_j=k$, $k_i=0$ for $i\neq j$ and $n=jk_j=jk$ with $j,k\in\NN$. Thus
$$\exists\;C,h,p'>0\;\forall\;j,k\in\NN:\;\;\;\frac{(jk)!}{(\lambda-f_{p}(0))^{k+1}}\left(\frac{R^{(p)}_j}{{\Gamma((2-\alpha)(j+1)+2)}}\right)^{k}\le Ch^{jk}M^{(p',\alpha)}_{jk}$$
is valid and clearly $(\lambda-f_{p}(0))^{k+1}\le h_1^{jk+1}$ for some $h_1>0$ (large) and all $k\in\NN_0$. Hence
\begin{equation}\label{theorem1siddiqiequ5}
	 \exists\;C,h,h_1,p'>0\;\forall\;j,k\in\NN:\;\;\;\left(\frac{R^{(p)}_j}{{\Gamma((2-\alpha)(j+1)+2)}}\right)^{k}\le Ch_1(hh_1)^{jk}\frac{M^{(p',\alpha)}_{jk}}{(jk)!}.
\end{equation}
By involving \eqref{theorem1siddiqiequ4} we estimate the left-hand side of \eqref{theorem1siddiqiequ5} as follows:
\begin{align*}
\frac{R^{(p)}_j}{{\Gamma((2-\alpha)(j+1)+2)}}& \ge\frac{j^{(1-\alpha)j}M^{(p,\alpha)}_j}{2^j{\Gamma((2-\alpha)(j+1)+2)}}\\ &\ge\frac{j!^{1-\alpha}M^{(p,\alpha)}_j}{2^j((2-\alpha)(j+1)+1){\Gamma((2-\alpha)(j+1)+1)}} \ge\frac{M^{(p,\alpha)}_j}{C_1{12}^jh_3^{j+1}j!}.
\end{align*}
The last estimate is valid since $(2-\alpha)(j+1)+1\le2(j+1)+(j+1)=3(j+1)\le{6}^j$ for all $j\in\NN$, and ${\Gamma((2-\alpha)(j+1)+1)}\le C_1h_2^{(2-\alpha)(j+1)}j!^{2-\alpha}$ for some $C_1,h_2\ge 1$ and all $j\ge 1$ (by the properties of the Gamma function), where we have put $h_3:=h_2^{2-\alpha}$. Consequently, by \eqref{theorem1siddiqiequ5} we get
\begin{equation*}%\label{theorem1siddiqiequ6}
	 \exists\;C,C_1,h,h_1,h_3,p'>0\;\forall\;j,k\in\NN:\;\;\;\left(\frac{M^{(p,\alpha)}_j}{j!}\right)^k\le Ch_1({12}hC_1h_1h^2_3)^{jk}\frac{M^{(p',\alpha)}_{jk}}{(jk)!},
\end{equation*}
and so
\begin{equation}\label{theorem1siddiqiequ7}
	\exists\;H\ge 1\;\exists\;p'(\ge p)>0\;\forall\;j,k\in\NN:\;\;\;\left(\frac{M^{(p,\alpha)}_j}{j!}\right)^{1/j}\le H\left(\frac{M^{(p',\alpha)}_{jk}}{(jk)!}\right)^{1/(jk)}.
\end{equation}
\eqref{theorem1siddiqiequ7} establishes \hyperlink{R-rai}{$(\mathcal{M}_{\{\on{rai}\}})$} for indices $p$ and $p'$ for all choices $j,k\in\NN$ and so for all multiplies $n=jk$ of $j\in\NN$. For the remaining cases let now $n\ge 1$ such that $jk<n<j(k+1)$ for some $j,k\in\NN$. Then, by using \eqref{theorem1siddiqiequ7} (with appearing constant $H$), \eqref{Stirling} and the fact that $j\mapsto(j^{(1-\alpha)j}M^{(p',\alpha)}_j)^{1/j}$
%$j\mapsto{(L^{(p')}_j)^{1/j}=(j^{(1-\alpha)j}M^{(p',\alpha)}_j)^{1/j}}$
is non-decreasing for each index $p'>0$ (by log-convexity), we estimate as follows:
\begin{align*}
	 \left(\frac{M^{(p',\alpha)}_{n}}{n!}\right)^{1/n}& =\frac{(n^{(1-\alpha)n}M^{(p',\alpha)}_{n})^{1/n}}{n^{1-\alpha}(n!)^{1/n}} \ge%{\color{red}\frac{(M^{(p')}_{jk})^{1/(jk)}}{n^{{1}-\alpha}(n!)^{1/n}}}=
	 \frac{((jk)^{(1-\alpha)jk}M^{(p',\alpha)}_{jk})^{1/(jk)}}{n^{1-\alpha}(n!)^{1/n}}
	\\&
	 =\frac{(jk)^{1-\alpha}}{n!^{1/n}n^{1-\alpha}} \left(\frac{M^{(p',\alpha)}_{jk}}{(jk)!}\right)^{1/(jk)}(jk)!^{1/(jk)}\\
&\ge\frac{1}{H}\left(\frac{M^{(p,\alpha)}_j}{j!}\right)^{1/j} \frac{(jk)!^{1/(jk)}}{n!^{1/n}}\left(\frac{jk}{n}\right)^{1-\alpha}
	\\&
	 \ge\frac{1}{H}\left(\frac{M^{(p,\alpha)}_j}{j!}\right)^{1/j} \frac{e^{-1}jk}{n}\left(\frac{jk}{j(k+1)}\right)^{1-\alpha}\\
&\ge\frac{1}{H}\left(\frac{M^{(p,\alpha)}_j}{j!}\right)^{1/j} \frac{jk}{ej(k+1)}\left(\frac{1}{2}\right)^{1-\alpha}	 \ge\frac{1}{He2^{2-\alpha}}\left(\frac{M^{(p,\alpha)}_j}{j!}\right)^{1/j}.
\end{align*}
Summarizing, property \hyperlink{R-rai}{$(\mathcal{M}_{\{\on{rai}\}})$} is verified for the matrix $\mathcal{M}^{\alpha}$ between the indices $p$ and $p'$ and when choosing the constant $C:=He2^{2-\alpha}(>H)$.\vspace{6pt}

$(a)\Rightarrow(e)$ This follows by $(ii)$ in Lemma \ref{condcomparison}.\vspace{6pt}

$(e)\Rightarrow(d)$ This follows by repeating the arguments in the proof of $(a)\Rightarrow(b)$ above (a word-by-word repetition of the proof in the ultradifferentiable setting), see \cite[Thm. 8.3.1]{dissertation}.\vspace{6pt}

$(d)\Rightarrow(b)$
%\red{Let $S$ be a sector in $\mathcal{R}$, property \hyperlink{R-Comega}{$(\mathcal{M}_{\{\text{C}^{\omega}\}})$} of $\mathcal{M}$ implies that the class $\mathcal{A}_{\{G^1\}}(S)$ is contained in $\mathcal{A}_{\{\mathcal{M}\}}(S)$ for any sector $S$.}
For all open set $U\subseteq \CC$, the property \hyperlink{R-Comega}{$(\mathcal{M}_{\{\text{C}^{\omega}\}})$} of $\mathcal{M}$ implies that the class $\mathcal{H}(U)$ is contained in $\mathcal{H}_{\{\mathcal{M}\}}(U)$. Since the class $\mathcal{A}_{\{\mathcal{M}\}}(S_{\alpha})$ is closed under composition, it is holomorphically closed too.
\qed\enddemo

\begin{remark}\label{ramimatrixdc}
\begin{itemize}
\item[(i)] If $\mathcal{M}$ has \hyperlink{R-dc}{$(\mathcal{M}_{\{\on{dc}\}})$} then $\mathcal{M}^\alpha$ has it too (the converse is not clear in general).
\item[(ii)] The condition that $\lim_{j\rightarrow+\infty}(j^{(1-\alpha)j}M^{(p)}_j)^{1/j}=\infty$ for all $p>0$ can be weakened as long as the log-convex regularization of $\overline{\G}^{1-\alpha}\M^{(p)}$ makes sense (for example, in case $\M^{(p)}=\overline{\G}^{\alpha-1}$). In this situation, the proof of Theorem~\ref{classequality} is still valid, Theorem~\ref{mainTthm2} can be applied and the availability of characteristic functions (needed in the previous proof of the implication $(c)\implies(a)$) is guaranteed. A similar comment can be made regarding the next corollary.
\end{itemize}
%	
%	If $\alpha=0$, we obtain the same result by replacing the sector by the positive real line in the theorems \ref{classequality} and \ref{theorem1siddiqi}.
\end{remark}

For a sequence $\M\in\RR_{>0}^{\NN_{0}}$ such that $\lim_{j\rightarrow+\infty}(j^{(1-\alpha)j}M^{(p)}_j)^{1/j}=\infty$, we can extend \cite[Thm. 1]{IderSiddiqi} by considering the constant weight matrix $\mathcal{M}=\{\M^{(p)}=\M: p>0\}$ and applying to it the previous result.

\begin{corollary}\label{coroCarlemanClassNarrowSector}
Let $\M\in\RR_{>0}^{\NN_{0}}$ be a sequence, and $0<\alpha\le 1$ be given such that $\lim_{j\rightarrow+\infty}(j^{(1-\alpha)j}M_j)^{1/j}=\infty$. Let $\M^{(\alpha)}:= \overline{\G}^{\alpha-1} \left(\overline{\G}^{1-\alpha}\M\right)^{\on{lc}}$. Then the following assertions are equivalent:
		\begin{itemize}
			\item[$(a)$] The sequence $\M^{(\alpha)}$ has the property \hyperlink{rai}{$(\on{rai})$}.
			
			\item[$(b)$] The class $\mathcal{A}_{\{\M\}}(S_{\alpha})$ is holomorphically closed.

			\item[$(c)$] The class $\mathcal{A}_{\{\M
				\}}(S_{\alpha})$ is inverse-closed.
			
		\end{itemize}
		If $\liminf_{j\rightarrow\infty}(\check{M}_j)^{1/j}>0$ and the sequence $\M^{(\alpha)}$ is (dc), then the list of equivalences can be extended by
		
		\begin{itemize}
			\item[$(d)$] The class $\mathcal{A}_{\{\M\}}(S_{\alpha})$ is closed under composition.

			\item[$(e)$] The sequence $\M^{(\alpha)}$ has the property \hyperlink{FdB}{$(\on{FdB})$}.
		\end{itemize}
\end{corollary}

\begin{remark}
    We may think of the situation for the ultradifferentiable class $\mathcal{E}_{\{\M\}}(0,+\infty)$, consisting of those complex-valued smooth functions on the half-line $(0,+\infty)$ subject to similar growth restrictions for their derivatives as in the ultraholomorphic case, as the limiting case when taking $\alpha=0$ in the previous result, i.e. when the sector $S_{\alpha}$ ``collapses'' to the ray $(0,+\infty)$. Then, it turns out that we (partially) recover the main result \cite[Thm.~1]{characterizationstabilitypaper}, see also \cite[Thm.~3.2]{compositionpaper}.
\end{remark}

Thanks to the construction of characteristic functions in classes defined in sectors of arbitrary opening, undertaken in Subsection~\ref{ssect.ConstrCharactFunct}, we study now the stability properties for {classes defined in sectors wider than a half-plane}.

\begin{theorem}\label{theorem1siddiqi.}
	Let $\mathcal{M}=\{\M^{(p)}: p>0\}$ be a weight matrix and consider $\alpha> 1$. For each $p>0$, we suppose that there exists some $\alpha_p>\alpha$ such that $\overline{\G}^{1-\alpha_p}\M^{(p)}$ is equivalent to a (lc) sequence $\LL^{(p)}$ depending on $\alpha_p$. Then the following assertions are equivalent:
	\begin{itemize}
		\item[$(a)$] The matrix $\mathcal{M}$ satisfies the property \hyperlink{R-rai}{$(\mathcal{M}_{\{\on{rai}\}})$}.
		
		\item[$(b)$] The class $\mathcal{A}_{\{\mathcal{M}\}}(S_{\alpha})$ is holomorphically closed.

		\item[$(c)$] The class $\mathcal{A}_{\{\mathcal{M}
			\}}(S_{\alpha})$ is inverse-closed.
		
	\end{itemize}
	If $\mathcal{M}$ has in addition \hyperlink{R-Comega}{$(\mathcal{M}_{\{\text{C}^{\omega}\}})$} and \hyperlink{R-dc}{$(\mathcal{M}_{\{\on{dc}\}})$}, then the list of equivalences can be extended by
	
	\begin{itemize}
		\item[$(d)$] The class $\mathcal{A}_{\{\mathcal{M}\}}(S_{\alpha})$ is closed under composition.
		
		\item[$(e)$] The matrix $\mathcal{M}$ satisfies the property \hyperlink{R-FdB}{$(\mathcal{M}_{\{\on{FdB}\}})$}.
	\end{itemize}
\end{theorem}

%Concerning the additional part we point out:
%
%
\demo{Proof} The proof of $(a)\Rightarrow(b)\Rightarrow(c)$ is similar to the one in Theorem~\ref{theorem1siddiqi}.

$(c)\Rightarrow(a)$ {Although the arguments are similar to those developed in the same implication in Theorem~\ref{theorem1siddiqi}, we consider it worthy to complete the details because now we will work with the original weight matrix (instead of $\mathcal{M}^\alpha$), and the characteristic functions are different in this framework.
}
%We follow the ideas from \cite[Theorem 1]{IderSiddiqi} and apply the constructions from the previous Section.
Let $p>0$ be arbitrary but from now on fixed. There exist  $\alpha_p>\alpha$ and $ \LL^{(p)} $ log-convex such that $\overline{\G}^{1-\alpha_p}\M^{(p)}\hyperlink{approx}{\approx} \LL^{(p)}$. Then, there exist $A_p,B_p>0$ such that $A_p^n n^{(1-\alpha_p)n}M^{(p)}_n\le L^{(p)}_n \le B_p^n n^{(1-\alpha_p)n}M^{(p)}_n  $ for all $n\in\NN_{0}$. According to Theorem \ref{mainTthm2} we put
\begin{equation*}%\label{theorem1siddiqiequ.}
	f_{p}(z):=\mathcal{T}_{\LL^{(p)}} (g_{\alpha,\alpha_p})(z).%\red{,\quad\alpha>1.}
\end{equation*}
By using \eqref{galphagoodestimate}, Lemma \ref{Scomplemma} and the above inequality we have
\begin{align*}
	|f^{(n)}_{p}(z)|&\le\sum_{k= 0}^\infty\frac{1}{2^k}L^{(p)}_k\frac{(\ell^{(p)}_{k})^n}{(\ell^{(p)}_{k})^k}|g_{\alpha,\alpha_p}^{(n)}(\ell^{(p)}_{k}z)|\le 2CD^nL^{(p)}_n\Gamma((\alpha_p-1)n+1)\\&
	\le E\widetilde{B}_p^nn^{(1-\alpha_p)n}M^{(p)}_nn^{(\alpha_p-1)n}= E\widetilde{B}_p^nM^{(p)}_n,
\end{align*}
for suitable constant $\widetilde{B}_p,C,D,E>1$ and for all $n\in\NN_0$ and $z\in S_{\alpha}$.
This estimate shows that $f_{p}\in\mathcal{A}_{\{\mathcal{M}\}}(S_{\alpha})$ and, in particular, it yields $\sup_{z\in S_{\alpha}}|f_{p}(z)|\leq E<+\infty$.

Set $R^{(p)}_n:=\sum_{k= 0}^\infty\frac{1}{2^k}L^{(p)}_k(\ell^{(p)}_{k})^{n-k}$, so that
\begin{equation}\label{theorem1siddiqiequ3.}
	\forall\;n\in\NN_0:\;\;\;f^{(n)}_{p}(0)=(-1)^n  \Gamma((\alpha_p-1)n+1)R_n^{(p)},
\end{equation}
and from Lemma \ref{Scomplemma},
\begin{equation}\label{theorem1siddiqiequ4.}
	 \forall\;n\in\NN_0:\;\;\;R^{(p)}_n\geq\frac{L^{(p)}_n}{2^n}\ge \frac{A_p^nn^{(1-\alpha_p)n}M^{(p)}_n}{2^n}.
\end{equation}
Now take $\lambda>E$ and put $\widetilde{f}_{p}:=\lambda-f_{p}$. Thus we get $\widetilde{f}_{p}\in\mathcal{A}_{\{\mathcal{M}\}}(S_{\alpha})$, and moreover $\inf_{z\in S_\alpha}|\widetilde{f}_{p}(z)|>0$. Since $\mathcal{A}_{\{\mathcal{M}\}}(S_{\alpha})$ is assumed to be inverse-closed, we get that $z\mapsto%\frac{1}{\widetilde{f}_{p}(z)}=
\frac{1}{\lambda-f_{p}(z)}\in\mathcal{A}_{\{\mathcal{M}\}}(S_{\alpha})$. When writing $g_{p} :z\mapsto\frac{1}{\lambda-z}$, the dependence on $p$ is justified because $\lambda$ is clearly depending on this chosen index. By applying the Fa\`a-di-Bruno-formula to the composition $g_{p} \circ f_{p}$ we get that for some $F,h>0$ and some index $p'>0$ (large) and for all $n\in\NN_0$,
\begin{align*}
	|(g_{p} \circ f_{p})^{(n)}(0)|&=\left|\sum_{\sum_{i=1}^nk_i=k, \sum_{i=1}^nik_i=n}\frac{n!}{k_1!\cdots k_n!} \frac{k!}{(\lambda-f_{p}(0))^{k+1}} \prod_{i=1}^n\left(\frac{f_p^{(i)}(0)}{i!}\right)^{k_i}\right|\\
&\le Fh^nM^{(p')}_n.
\end{align*}
Using \eqref{theorem1siddiqiequ3.}
%we see $\left(\frac{f_{p}^{(i)}(0)}{i!}\right)^{k_i}= \left(\frac{(-1)^i\Gamma((\alpha_p-1)i+1)S_i^{p}}{i!}\right)^{k_i}$ for all $1\le i\le n$,
and since $ \prod_{i=1}^n (-1)^{ik_i}=(-1)^n $, we deduce that for every $n\in\NN_0$
\begin{align*}
\sum_{\sum_{i=1}^nk_i=k, \sum_{i=1}^nik_i=n} \frac{n!}{k_1!\cdots k_n!} \frac{k!}{(\lambda-f_{p}(0))^{k+1}} \prod_{i=1}^n\left(\frac{\Gamma((\alpha_p-1)i+1)R_i^{(p)}}{i!}\right)^{k_i}
\le Fh^nM^{(p')}_n.
\end{align*}
Given $j,k\in\NN$, we focus on the summand for $k_j=k$, $k_i=0$ for $i\neq j$ and $n=jk_j=jk$, so we get that
$$\exists\;F,h,p'>0\;\forall\;j,k\in\NN:\;\;\;\frac{(jk)!}{(\lambda-f_{p}(0))^{k+1}} \left(\frac{\Gamma((\alpha_p-1)j+1)R_j^{(p)}}{j!}\right)^{k}\le Fh^{jk}M^{(p')}_{jk}.
$$
Clearly, $(\lambda-f_{p}(0))^{k+1}\le h_1^{jk+1}$ for some $h_1>0$ (large) and all $k\in\NN_0$. Hence, for all $j,k\in\NN$ we have
\begin{equation}\label{theorem1siddiqiequ5.}
\exists\;F,h,h_1,p'>0\;\forall\;j,k\in\NN:\;\;\;\left(\frac{\Gamma((\alpha_p-1)j+1)R_j^{(p)}}{j!}\right)^{k}\le Fh_1(hh_1)^{jk}\frac{M^{(p')}_{jk}}{(jk)!}.
\end{equation}
By involving \eqref{theorem1siddiqiequ4.} we estimate the left-hand side of \eqref{theorem1siddiqiequ5.} as follows:
\begin{align*}
\frac{\Gamma((\alpha_p-1)j+1)R_j^{(p)}}{j!}& \ge\frac{A_p^jj^{(1-\alpha_p)j}\Gamma((\alpha_p-1)j+1)M^{(p)}_j}{2^jj!}
	\\&
\ge\frac{{\widetilde{A}^j_p}j^{(1-\alpha_p)j}j^{(\alpha_p-1)j}M^{(p)}_j}{2^jj!}= \frac{M^{(p)}_j}{{\overline{A}^j_p}j!}.
\end{align*}
The last inequality is a consequence of the properties of the Gamma function for a suitable constant $\widetilde{A}_p>0$, and we have put $\overline{A}_p=2/\widetilde{A}_p$. Consequently, by \eqref{theorem1siddiqiequ5.} we get
\begin{equation*}%\label{theorem1siddiqiequ6.}
\exists\;F,h,h_1,\overline{A}_p,p'>0\;\forall\;j,k\in\NN:\;\;\;\left(\frac{M^{(p)}_j}{j!}\right)^k\le Fh_1(hh_1\overline{A}_p)^{jk}\frac{M^{(p')}_{jk}}{(jk)!},
\end{equation*}
and so there exists $H\ge 1$ such that
\begin{equation}\label{theorem1siddiqiequ7.}
\left(\frac{M^{(p)}_j}{j!}\right)^{1/j}\le H\left(\frac{M^{(p')}_{jk}}{(jk)!}\right)^{1/(jk)}.
\end{equation}
{Equation} \eqref{theorem1siddiqiequ7.} establishes \hyperlink{R-rai}{$(\mathcal{M}_{\{\on{rai}\}})$} for indices $p$ and $p'$ for all choices $j,k\in\NN$ and so for all multiples $n=jk$ of $j\in\NN$. For the remaining cases let now $n\ge 1$ such that $jk<n<j(k+1)$ for some $j,k\in\NN$. Then, by using \eqref{theorem1siddiqiequ7.}, \eqref{Stirling}, the equivalence $\overline{\G}^{1-\alpha_{p'}}\M^{(p')}\hyperlink{approx}{\approx} \LL^{(p')}$ and the fact that $j\mapsto(L_j^{(p')})^{1/j}$
is non-decreasing for each index $p'>0$, we estimate
\begin{align*}
\left(\frac{M^{(p')}_{n}}{n!}\right)^{1/n}& =\frac{(B_{p'}^nn^{(1-\alpha_{p'})n}M^{(p')}_{n})^{1/n}}{B_{p'}n^{1-\alpha_{p'}}(n!)^{1/n}} \ge\frac{(L^{(p')}_{n})^{1/n}}{B_{p'}n^{1-\alpha_{p'}}(n!)^{1/n}} \ge\frac{(L^{(p')}_{jk})^{1/(jk)}}{B_{p'}n^{1-\alpha_{p'}}(n!)^{1/n}}
	\\& \ge%{\color{red}\frac{(M^{(p')}_{jk})^{1/(jk)}}{n^{{1}-\alpha}(n!)^{1/n}}}=
\frac{(A_{p'}^{jk}(jk)^{(1-\alpha_{p'})jk}M^{(p')}_{jk})^{1/(jk)}}{B_{p'}n^{1-\alpha_{p'}}(n!)^{1/n}}
=\frac{A_{p'}(jk)^{1-\alpha_{p'}}}{B_{p'}n!^{1/n}n^{1-\alpha_{p'}}} \left(\frac{M^{(p')}_{jk}}{(jk)!}\right)^{1/(jk)}(jk)!^{1/(jk)}
\\&
\ge\frac{A_{p'}}{B_{p'}H}\left(\frac{M^{(p)}_j}{j!}\right)^{1/j} \frac{(jk)!^{1/(jk)}}{n!^{1/n}}\left(\frac{jk}{n}\right)^{1-\alpha_{p'}}
\\&
\ge\frac{A_{p'}}{B_{p'}H}\left(\frac{M^{(p)}_j}{j!}\right)^{1/j} \frac{e^{-1}jk}{n}%\red{\left(\frac{jk}{j(k+1)}\right)^{1-\alpha_{p'}}}
	\\&
\ge\frac{A_{p'}}{B_{p'}H}\left(\frac{M^{(p)}_j}{j!}\right)^{1/j} \frac{jk}{ej(k+1)}%\red{\left(\frac{1}{2}\right)^{1-\alpha_{p'}}}
\ge\frac{A_{p'}}{2B_{p'}He} \left(\frac{M^{(p)}_j}{j!}\right)^{1/j}.
\end{align*}
Summarizing, property \hyperlink{R-rai}{$(\mathcal{M}_{\{\on{rai}\}})$} is verified for the matrix $\mathcal{M}$ between the indices $p$ and $p'$ and when choosing the constant $C:=2B_{p'}He/A_{p'}$.\vspace{6pt}

$(a)\Rightarrow(e)$ and $(d)\Rightarrow(b)$ are as in Theorem~\ref{theorem1siddiqi}.
%This follows by $(ii)$ in Lemma \ref{condcomparison}.
\vspace{6pt}

$(e)\Rightarrow(d)$ One can repeat the proof in the ultradifferentiable setting, see~\cite[Thm. 8.3.1]{dissertation}.
\qed\enddemo
%noting that
% \overline{\G}^{1-\alpha_{p}}\M^{(p)}\hyperlink{approx}{\approx} \LL^{(p)} $ implies that
%$\M^{(p)}\hyperlink{approx}{\approx} \overline{\G}^{\alpha_{p}-1}\LL^{(p)}$, and so $\M^{(p)}$ is equivalent to an (lc) sequence, see \cite[Thm. 4.9 $(4)\Rightarrow(1)$]{compositionpaper}.

%%$(d)\Rightarrow(b)$ \red{Let $S$ be a sector in $\mathcal{R}$, property \hyperlink{R-Comega}{$(\mathcal{M}_{\{\text{C}^{\omega}\}})$} of $\mathcal{M}$ implies that the class $\mathcal{A}_{\{G^1\}}(S)$ is contained in $\mathcal{A}_{\{\mathcal{M}\}}(S)$ for any sector $S$.}
%%{For all open set $U\subseteq \CC$, the property \hyperlink{R-Comega}{$(\mathcal{M}_{\{\text{C}^{\omega}\}})$} of $\mathcal{M}$ implies that the class $\mathcal{H}(U)$ is contained in $\mathcal{H}_{\{\mathcal{M}\}}(U)$}. Finally, the previous comment and the fact that the class $\mathcal{A}_{\{\mathcal{M}\}}(S_{\alpha})$ is closed under composition implies that the class is holomorphically closed too.

%Let $S$ be a sector in $\mathcal{R}$, property \hyperlink{R-Comega}{$(\mathcal{M}_{\{\text{C}^{\omega}\}})$} of $\mathcal{M}$ implies that the class $\mathcal{A}_{\{G^1\}}(S)$ is contained in $\mathcal{A}_{\{\mathcal{M}\}}(S)$ for any sector $S$. Finally, the previous comment and the fact that the class $\mathcal{A}_{\{\mathcal{M}\}}(S_{\alpha})$ is closed under composition implies that the class is holomorphically closed too.

\begin{remark}\label{ramimatrixdc.}
In the same line of Remark~\ref{rem.gammaindexcharactfunct}, if for a weight matrix $\mathcal{M}=\{\M^{(p)}: p>0\}$  we know that $\gamma(\M^{(p)})>\alpha-1$ for all $p>0$, then there exist{s} some $\alpha_p>\alpha$ such that $\overline{\G}^{1-\alpha_p}\M^{(p)}$ is equivalent to a (lc) sequence $\LL^{(p)}$ depending on $\alpha_p$.
\end{remark}

%\begin{example}
%{We can also add some examples of weight matrix for which we can apply this result, e.g. let $\alpha>1$ and we consider the family of weight matrix $$\mathcal{M}_\beta=\{\M^{(p)}=\overline{\G}^{\beta+p}: p>0\}\qquad \forall\beta>\alpha-1.$$}
%\end{example}

Note that there exist some differences between the statements of the Theorems \ref{theorem1siddiqi} and \ref{theorem1siddiqi.}, concerning the fact that the conditions for stability are imposed on different weight matrices, $\mathcal{M}$ or $\mathcal{M^\alpha}$. In general, if $\alpha>1$ we only know that $\mathcal{A}_{\{\mathcal{M^\alpha}\}}(S_{\alpha})\subset \mathcal{A}_{\{\mathcal{M}\}}(S_{\alpha})$. However, {the hypotheses of the second theorem have strong implications and, under an additional assumption, these results perfectly match, as the next proposition shows.}\\
\begin{proposition} Let $\mathcal{M}=\{\M^{(p)}: p>0\}$ be a given weight matrix. Suppose that {for every $p>0$ there exists  $\alpha_p>0$ such that $\overline{\G}^{1-\alpha_p}\M^{(p)}$ is equivalent to a (lc) sequence $\LL^{(p)}$, and that there exists $\beta\in\RR$ such that $\beta<\alpha_p$ for all $p>0$. Then, for every $p>0$ one has $\lim_{j\rightarrow+\infty}(j^{(1-\beta)j}M^{(p)}_j)^{1/j}=\infty$, $\mathcal{M}$ and $\mathcal{M}^{\beta}$ (defined as in \eqref{ramiequ}) are $R$-equivalent, and therefore $\mathcal{M}$ satisfies the property \hyperlink{R-rai}{$(\mathcal{M}_{\{\on{rai}\}})$} (resp.\hyperlink{R-FdB}{$(\mathcal{M}_{\{\on{FdB}\}})$}) if and only if the matrix $\mathcal{M}^{\beta}$ satisfies this condition too. Moreover, $\mathcal{A}_{\{\mathcal{M^{\text{$\beta$}}}\}}(S_{\gamma})= \mathcal{A}_{\{\mathcal{M}\}}(S_{\gamma})$, for all $\gamma>0$.}
\end{proposition}
\demo{Proof} {Let $p>0$ be arbitrary but fixed. First, note that $$\overline{\G}^{1-\beta}\M^{(p)}= \overline{\G}^{\alpha_p-\beta}(\overline{\G}^{1-\alpha_p}\M^{(p)}) \hyperlink{approx}{\approx}\overline{\G}^{\alpha_p-\beta}\LL^{(p)}=:\widetilde{\LL}^{(p)},$$
where the sequence $\widetilde{\LL}^{(p)}$ is log-convex (as the product of two such sequences).

On the one hand, the condition $\LL^{(p)}\hyperlink{approx}{\approx}\overline{\G}^{1-\alpha_p}\M^{(p)}$ guarantees that there exists some $A>0$ such that $A^jL_j^{(p)}\leq j^{(1-\alpha_p)j}M^{(p)}_j$, for all $j\in\NN_0$. Moreover, let us observe that for all $j>0$, we can estimate $(j^{(1-\beta)j}M^{(p)}_j)^{1/j}=j^{(\alpha_p-\beta)}(j^{(1-\alpha_p)j}M^{(p)}_j)^{1/j}\geq j^{(\alpha_p-\beta)}(A^jL^{(p)}_j)^{1/j}$, and thanks to the fact that $\LL^{(p)}$ is (lc) and $\alpha_p>\beta$, we deduce that $\lim_{j\rightarrow+\infty}(j^{(1-\beta)j}M^{(p)}_j)^{1/j}=\infty$. Moreover, there exists some $\widetilde{A}>0$ such that the (lc) sequence $\mathbb{B}^{(p)}:=(\widetilde{A}^j\widetilde{L}_j^{(p)})_j$ satisfies $\mathbb{B}^{(p)}\leq \overline{\G}^{1-\beta}\M^{(p)}$.
Then, we have that $\mathbb{B}^{(p)}=(\mathbb{B}^{(p)})^{\on{lc}}\leq (\overline{\G}^{1-\beta}\M^{(p)})^{\on{lc}}$, which implies that $\widetilde{\LL}^{(p)}\hyperlink{preceq}{\preceq}(\overline{\G}^{1-\beta}\M^{(p)})^{\on{lc}}$.

On the other hand, we observe that $\overline{\G}^{1-\beta}\M^{(p)}\hyperlink{preceq}{\preceq}\widetilde{\LL}^{(p)}$, and therefore,  $(\overline{\G}^{1-\beta}\M^{(p)})^{\on{lc}}\hyperlink{preceq}{\preceq}\widetilde{\LL}^{(p)}$. Finally, we conclude that $\widetilde{\LL}^{(p)}\hyperlink{approx}{\approx}(\overline{\G}^{1-\beta}\M^{(p)})^{\on{lc}}$.

The previous equivalence ensures that $\M^{(p,\beta)}$ is equivalent to $\overline{\G}^{\beta-1}\widetilde{\LL}^{(p)}$, and therefore $\M^{(p)}\hyperlink{approx}{\approx}\M^{(p,\beta)}$. %
%, and in particular by taking $\beta=\alpha$ we have that $\M^{(p)}\hyperlink{approx}{\approx}\M^{(p,\alpha)}$.
Finally, the two matrices $\mathcal{M}$ and $\mathcal{M}^{\beta}$ are $R$-equivalent, and the property \hyperlink{R-rai}{$(\mathcal{M}_{\{\on{rai}\}})$} (resp.\hyperlink{R-FdB}{$(\mathcal{M}_{\{\on{FdB}\}})$})
is stable under $R$-equivalence, see \cite[Remark 8.2.2]{dissertation}.}
\qed\enddemo

{Under the assumptions of the previous proposition, we can prove a weaker variant of Theorem \ref{theorem1siddiqi.} using a similar technique to the one used in the proof of Theorem \ref{theorem1siddiqi}.}
%{Option 2: Thanks to the construction of characteristic functions in classes defined in sectors of arbitrary opening, undertaken in Subsection~\ref{ssect.ConstrCharactFunct}, we study now the stability properties for such classes. Or do we erase it?}
%Thanks to the previous Proposition, we can adapt the proof for $(c)\Rightarrow(a)$. We include the proof, but we should omit them in the final version
\begin{corollary}
Let $\mathcal{M}=\{\M^{(p)}: p>0\}$ be a weight matrix and consider $\alpha> 1$. For each $p>0$, we suppose that there exist some $\alpha_p>\alpha$ such that $\overline{\G}^{1-\alpha_p}\M^{(p)}$ is equivalent to a (lc) sequence $\LL^{(p)}$ depending on $\alpha_p$, {and that there exists $\beta>\alpha$ such that $\beta<\alpha_p$ for all $p>0$}. Then the following assertions are equivalent:
	\begin{itemize}
		\item[$(a)$] The matrix $\mathcal{M}$, {or equivalently $\mathcal{M}^\beta$}, satisfies property \hyperlink{R-rai}{$(\mathcal{M}_{\{\on{rai}\}})$}.
		
		\item[$(b)$] The class $\mathcal{A}_{\{\mathcal{M}\}}(S_{\alpha})$ is holomorphically closed.

		\item[$(c)$] The class $\mathcal{A}_{\{\mathcal{M}
			\}}(S_{\alpha})$ is inverse-closed.
		
	\end{itemize}
	If $\mathcal{M}$ has in addition \hyperlink{R-Comega}{$(\mathcal{M}_{\{\text{C}^{\omega}\}})$} and \hyperlink{R-dc}{$(\mathcal{M}_{\{\on{dc}\}})$}, then the list of equivalences can be extended by
	
	\begin{itemize}
		\item[$(d)$] The class $\mathcal{A}_{\{\mathcal{M}\}}(S_{\alpha})$ is closed under composition.
		
		\item[$(e)$] The matrix $\mathcal{M}$, {or equivalently $\mathcal{M}^\beta$}, satisfies  property \hyperlink{R-FdB}{$(\mathcal{M}_{\{\on{FdB}\}})$}.
	\end{itemize}
\end{corollary}

We end this section by providing the version of Corollary~\ref{coroCarlemanClassNarrowSector} for wide sectors, which can be again deduced as a straightforward consequence of the corresponding result for weight matrices, Theorem~\ref{theorem1siddiqi.}.

\begin{corollary}\label{coroCarlemanClassWideSector} Let $\M\in\RR_{>0}^{\NN_{0}}$ and $\alpha> 1$. Suppose there exists $\alpha'>\alpha$ such that $\overline{\G}^{1-\alpha'}\M$ is equivalent to an (lc) sequence $\LL$ (depending on $\alpha'$).
%$\LL^{(\alpha)}=\overline{\G}^{\alpha-1} \left(\overline{\G}^{1-\alpha}\LL\right)^{\on{lc}}$.
Then the following assertions are equivalent:
		\begin{itemize}
			\item[$(a)$] The sequence $\M$ has the property \hyperlink{rai}{$(\on{rai})$}.
			
			\item[$(b)$] The class $\mathcal{A}_{\{\M\}}(S_{\alpha})$ is holomorphically closed.

			\item[$(c)$] The class $\mathcal{A}_{\{\M
				\}}(S_{\alpha})$ is inverse-closed.
			
		\end{itemize}
		If $\liminf_{j\rightarrow\infty}(\check{M}_j)^{1/j}>0$ and $\M$ is (dc), then the list of equivalences can be extended by
		
		\begin{itemize}
			\item[$(d)$] The class $\mathcal{A}_{\{\M\}}(S_{\alpha})$ is closed under composition.

			\item[$(e)$] The sequence $\M$ has the property \hyperlink{FdB}{$(\on{FdB})$}.
		\end{itemize}
\end{corollary}
\section{The weight function case}\label{sect.StabilityWeightFUnctCase}

We start proving, for the reader's convenience, how the condition \hyperlink{R-rai}{$(\mathcal{M}_{\{\on{rai}\}})$} for a weight matrix associated to a weight function $\omega$ translates into a condition on $\omega$. {Note that this matrix has \hyperlink{Mlc}{$(\mathcal{M}_{\on{lc}})$} and therefore $(\mathcal{M}_{\omega})^{\alpha}\equiv\mathcal{M}_{\omega}$ for all $\alpha\in(0,1]$.}

\begin{lemma}\label{rai-alpha0}
Let {$\omega\in\hyperlink{omset0}{\mathcal{W}_0}$} be given with associated weight matrix $\mathcal{M}_{\omega}:=\{\mathbb{W}^{(\ell)}: \ell>0\}$. Then the following are equivalent:
	\begin{itemize}
		\item[$(a)$] The matrix $\mathcal{M}_{\omega}$ has \hyperlink{R-rai}{$(\mathcal{M}_{\{\on{rai}\}})$}, i.e. (recall $\check{W}^{(\ell)}_j=W^{(\ell)}_j/j!$)
		$$\forall\;\ell>0\;\exists\;\ell'>0\;\exists\;H\ge 1\;\forall\;1\le j\le k:\;\;\;(\check{W}^{(\ell)}_j)^{1/j}\le H(\check{W}^{(\ell')}_k)^{1/k}.$$
		
		\item[$(b)$] $\omega$ satisfies the condition $(\alpha_0)$ (see \eqref{alpha0}), i.e.
		$$\exists\;C\ge 1\;\exists\;t_0\ge 0\;\forall\;\lambda\ge 1\;\forall\;t\ge t_0:\;\;\;\omega(\lambda t)\le C\lambda\omega(t).$$
  \end{itemize}
\end{lemma}
\demo{Proof}
$(a)\Rightarrow(b)$ The property \hyperlink{R-rai}{$(\mathcal{M}_{\{\on{rai}\}})$} is preserved under equivalence of matrices, then $\mathcal{M}_{\omega_{{\W}^{(\ell)}}}$ has \hyperlink{R-rai}{$(\mathcal{M}_{\{\on{rai}\}})$} for some/any $l>0$. By \cite[Thm. 4.5 $(iv)\Leftrightarrow(i)$]{subaddlike} $\omega_{{\W}^{(\ell)}}$ satisfies the condition $(\alpha_0)$, and therefore $\omega$ satisfies it too, because $\omega\sim\omega_{{\W}^{(\ell)}}$ (see~\eqref{goodequivalenceclassic}) and the condition $(\alpha_0)$ is preserved under equivalence of weight functions.

$(b)\Rightarrow(a)$ If $\omega$ satisfies the condition $(\alpha_0)$, then $\omega_{{\W}^{(\ell)}}$ satisfies it too (arguing as before). By \cite[Thm. 4.5 $(i)\Leftrightarrow(iv)$]{subaddlike}, the matrix $\mathcal{M}_{\omega_{{\W}^{(\ell)}}}$ has \hyperlink{R-rai}{$(\mathcal{M}_{\{\on{rai}\}})$} for some/any $l>0$. Finally, by \cite[Lemma 5.3.1]{dissertation} the matrices $\mathcal{M}_{\omega_{{\W}^{(\ell)}}}$ and $\mathcal{M}_{\omega} $ are equivalent, and \hyperlink{R-rai}{$(\mathcal{M}_{\{\on{rai}\}})$} is preserved under equivalence of matrices.
%}
\qed\enddemo
\vspace{6pt}
%\red{CASE $0<\alpha\leq 1$:}
%\vspace{6pt}

We can provide now a statement about stability properties for classes associated to a weight function in small sectors.

\begin{theorem}\label{theorem1siddiqcase1}
Let $\omega\in\hyperlink{omset1}{\mathcal{W}}$ be given with associated weight matrix $\mathcal{M}_{\omega}:=\{\mathbb{W}^{(\ell)}: \ell>0\}$ and let $0<\alpha\le 1$. Then the following are equivalent:
	\begin{itemize}
		\item[$(a)$] The matrix $\mathcal{M}_{\omega}$ has \hyperlink{R-rai}{$(\mathcal{M}_{\{\on{rai}\}})$}.
		
		\item[$(b)$] $\omega$ satisfies the condition $(\alpha_0)$ (see \eqref{alpha0}).
		
		\item[$(c)$] The class $\mathcal{A}_{\{\omega\}}(S_{\alpha})$ is holomorphically closed.
		
		\item[$(d)$] The class $\mathcal{A}_{\{\omega\}}(S_{\alpha})$ is inverse-closed.
  \end{itemize}
	If $\omega$ has in addition $(\omega_2)$, then the list of equivalences can be extended by:

 \begin{itemize}
		\item[$(e)$] The class $\mathcal{A}_{\{\omega\}}(S_{\alpha})$ is closed under composition.
		
		\item[$(f)$] The matrix $\mathcal{M}_{\omega}$ satisfies the property \hyperlink{R-FdB}{$(\mathcal{M}_{\{\on{FdB}\}})$}.
\end{itemize}	
\end{theorem}
\demo{Proof}
The equivalence $(a)\Leftrightarrow(b)$ is a consequence of the Lemma \ref{rai-alpha0}. Moreover, the equivalences $(a)\Leftrightarrow(c)\Leftrightarrow(d)\Leftrightarrow(e)\Leftrightarrow(f)$ follow by applying Theorem \ref{theorem1siddiqi} to $\mathcal{M}\equiv\mathcal{M}_{\omega}$. Let us observe that $\mathcal{M}^{\alpha}\equiv\mathcal{M}_{\omega}$, thanks to the fact that $\W^{(\ell)}$ is (lc) for all $\ell>0$. Moreover, $\omega$ has $(\omega_1)$ and therefore $\mathcal{A}_{\{\omega\}}(S_{\alpha})=\mathcal{A}_{\{\mathcal{M}_{\omega}\}}(S_{\alpha})$, see \eqref{equaEqualitySpacesWeightFunctionMatrix}. In addition, note that $\mathcal{M}_{\omega}$ has automatically \hyperlink{R-dc}{$(\mathcal{M}_{\{\on{dc}\}})$} by \eqref{newmoderategrowth}.
\qed\enddemo

\begin{remark}
When taking $\alpha=0$ in the previous result, i.e., when the sector $S_{\alpha}$ ''collapses'' to the ray $(0,+\infty)$, then we (partially) get back the main result \cite[Thm. 3]{characterizationstabilitypaper} for the ultradifferentiable class $\mathcal{E}_{\{\omega\}}((0,+\infty))$, see also \cite[Thm. 6.3]{compositionpaper}.
\end{remark}

%\red{EN ESTE CASO QUEREMOS QUE $\alpha>1$. EL TEOREMA SER\'IA ESTE:}\\
%{El comentario diciendo que si tienes la condición $(\alpha_0)$ para $\omega$ entonces tienes $(\omega_1)$ para ella está incluído en la demostración. Luego podemos reemplazar $\omega\in\hyperlink{omset1}{\mathcal{W}}$ por $\omega\in\hyperlink{omset0}{\mathcal{W}_0}$.}\vspace{6pt}

The next lemma will be necessary for stating a similar result for wide sectors.

\begin{lemma}\label{ramifalpha0}
Let $\omega\in\hyperlink{omset0}{\mathcal{W}_0}$ be given with associated weight matrix $\mathcal{M}_{\omega}:=\{\mathbb{W}^{(\ell)}: \ell>0\}$. Suppose there exists $s > 0$ such that, for $\omega^s(t):=\omega(t^s)$, one has:\begin{itemize}
    \item[(i)]  $\omega^s(t) = o(t)$ as $t\to\infty$, (i.e., $\omega^s(t)$ has $(\omega_5)$.)
    \item[(ii)]  $\omega^s$ satisfies the condition $(\alpha_0)$, i.e., it is
 equivalent to a concave weight function.
\end{itemize}
%
%Then $\mathcal{U}=\{\mathbb{U}^{(\ell)}:= \G^{s}[(\check{\mathbb{V}}^{(\ell,s)})^{\on{lc}}]^s: \ell>0\}$ are $R$-equivalent, and for each $\ell>0$, the sequence $\overline{\G}^{-s}\mathbb{U}^{(\ell)}$ is equivalent to a (lc) sequence $\LL^{(\ell)}$ depending on $s$.
Then there exists a weight matrix $\,\mathcal{U}=\{\mathbb{U}^{(\ell)}: \ell>0\}$, $R$-equivalent to $\mathcal{M}_{\omega}$, and such that for each $\ell>0$, the sequence $\overline{\G}^{\;-s}\mathbb{U}^{(\ell)}$ is equivalent to an (lc) sequence $\LL^{(\ell)}$ depending on~$s$.
\end{lemma}
\demo{Proof}
First, let us consider the matrix $\mathcal{M}_{\omega^s}:=\{\mathbb{V}^{(\ell,s)}: \ell>0\}$. There exists a relation between both matrices (see \cite{sectorialextensions1}), more precisely, for all $\ell>0$ we have that $\mathbb{V}^{(\ell,s)}=(\mathbb{W}^{(\ell/s)})^{1/s}$. So, we can write
\begin{equation*}
    \mathbb{W}^{(\ell)}=(\mathbb{V}^{(\ell s,s)})^s=\G^{s}(\check{\mathbb{V}}^{(\ell s,s)})^s \qquad \ell>0.
\end{equation*}
Now, by taking into account that $\omega^s$ satisfies the condition $(\alpha_0)$ and $(\omega_5)$ we deduce from \cite[Prop 3]{whitneyextensionmixedweightfunctionII} that the matrices $\check{\mathcal{M}}_{\omega^s}:=\{\check{\mathbb{V}}^{(\ell,s)}: \ell>0\}$ and  $\check{\mathcal{M}}_{\omega^s}^{\on{lc}}$ are $R$-equivalent.
Finally, since taking the power $s$ in each sequence of these two matrices respects $R$-equivalence for the resulting matrices, we deduce that  $\mathcal{U}:=\{\G^{s}[(\check{\mathbb{V}}^{(\ell,s)})^{\on{lc}}]^s: \ell>0\}$ and $\mathcal{M}_{\omega}$ are $R$-equivalent.
\qed\enddemo

\begin{theorem}\label{theoremWeightFunctionWideSector}
Let $\omega\in\hyperlink{omset0}{\mathcal{W}_0}$ be given with associated weight matrix $\mathcal{M}_{\omega}:=\{\mathbb{W}^{(\ell)}: \ell>0\}$ and let $\alpha> 1$. Suppose there exists $s > \alpha -1$ such that, for $\omega^s(t):=\omega(t^s)$, one has:\begin{itemize}
    \item[(i)]  $\omega^s(t) = o(t)$ as $t\to\infty$, (i.e $\omega^s(t)$ has $(\omega_5)$).
    \item[(ii)]  $\omega^s$ satisfies the condition $(\alpha_0)$, i.e., it is
 equivalent to a concave weight function.
\end{itemize}
Then the following are equivalent:
	\begin{itemize}
		\item[$(a)$] The matrix $\mathcal{M}_{\omega}$ has \hyperlink{R-rai}{$(\mathcal{M}_{\{\on{rai}\}})$}.

		\item[$(b)$] $\omega$ satisfies the condition $(\alpha_0)$.
  % (see \eqref{alpha0}), so $$\exists\;C\ge 1\;\exists\;t_0\ge 0\;\forall\;\lambda\ge 1\;\forall\;t\ge t_0:\;\;\;\omega(\lambda t)\le C\lambda\omega(t).$$
		
		\item[$(c)$] The class $\mathcal{A}_{\{\omega\}}(S_{\alpha})$ is holomorphically closed.
		
		\item[$(d)$] The class $\mathcal{A}_{\{\omega\}}(S_{\alpha})$ is inverse-closed.
  \end{itemize}
	If $\omega$ has in addition $(\omega_2)$, then the list of equivalences can be extended by:

 \begin{itemize}
		\item[$(e)$] The class $\mathcal{A}_{\{\omega\}}(S_{\alpha})$ is closed under composition.
		
		\item[$(f)$] The matrix $\mathcal{M}_{\omega}$ satisfies the condition \hyperlink{R-FdB}{$(\mathcal{M}_{\{\on{FdB}\}})$}.
\end{itemize}	
\end{theorem}

\demo{Proof}
The equivalence $(a)\Leftrightarrow(b)$ is a consequence of Lemma \ref{rai-alpha0}. Lemma \ref{ramifalpha0} ensures that there exists a weight matrix $\,\mathcal{U}:=\{\mathbb{U}^{(\ell)}: \ell>0\}$, $R$-equivalent to $\mathcal{M}_{\omega}$ (and therefore $\mathcal{A}_{\{\mathcal{U}\}}(S_{\alpha})= \mathcal{A}_{\{\mathcal{M}_{\omega}\}}(S_{\alpha})$), such that for each $\ell>0$ the sequence $\overline{\G}^{\;-s}\mathbb{U}^{(\ell)}$ is equivalent to a (lc) sequence $\LL^{(\ell)}$ depending on $s$. Then, the equivalences $(a)\Leftrightarrow(c)\Leftrightarrow(d)\Leftrightarrow(e) \Leftrightarrow(f)$ follow by applying Theorem \ref{theorem1siddiqi.} to $\mathcal{M}\equiv\mathcal{U}$, and taking $\alpha_\ell=s+1$. Finally, thanks to the fact that $\omega^s$ has $(\alpha_0)$, then $\omega$ satisfies $(\omega_1)$ and therefore $\mathcal{A}_{\{\omega\}}(S_{\alpha})= \mathcal{A}_{\{\mathcal{M}_{\omega}\}}(S_{\alpha})$, see \eqref{equaEqualitySpacesWeightFunctionMatrix}.
In addition, note that $\mathcal{M}_{\omega}$ has automatically \hyperlink{R-dc}{$(\mathcal{M}_{\{\on{dc}\}})$} by \eqref{newmoderategrowth}. And $(\omega_2)$ for $\omega$ implies that $\mathcal{M}_\omega$ has \hyperlink{holom}{$(\mathcal{M}_{\mathcal{H}})$}. Finally, the conditions \hyperlink{R-dc}{$(\mathcal{M}_{\{\on{dc}\}})$} and \hyperlink{holom}{$(\mathcal{M}_{\mathcal{H}})$} are stable under $R$-equivalence, and therefore $\mathcal{U}$ satisfies both too.
\qed\enddemo

\begin{remark}\label{gammaindex-omega}
The hypotheses (i) and (ii) on $\omega$ in Theorem~\ref{theoremWeightFunctionWideSector} can be quickly guaranteed by the condition $\gamma(\omega) > \alpha -1$, in terms of the index described in Subsection~\ref{ssect.IndexGammaOmega}. Note that, by choosing $s$ such that $\gamma(\omega) > s > \alpha-1$,
we have  $\gamma(\omega^s) = \gamma(\omega) / s > 1$ (see property $(iii)$ in that subsection), and this fact implies:
\begin{itemize}
    \item[(a)] By \cite[Remark 2.15 $(i)\Rightarrow(v)$]{index}, we have property $(\omega_5)$ for $\omega^s$.

    \item[(b)] By \cite[Thm. 2.11 $(v)\Rightarrow(ii)$]{index}, we deduce that $\omega^s$ is equivalent to a concave weight function, and so $(\alpha_0)$ is satisfied by $\omega^s$.
\end{itemize}
\end{remark}

\begin{remark}
	%\red{Under suitable conditions, we can prove the equivalence between all the statements in the previous theorem. More precisely:}
{In some situations it is straightforward that all the conditions on the weight function $\omega$ in the previous result are satisfied, and so all the statements (a) through (f) are equivalent. We comment on two special cases:}
\begin{itemize}
    \item[(i)] If $2>\alpha> 1$, suppose that $\omega(t) = O(t)$ as $t\to\infty$, (i.e $\omega(t)$ has $(\omega_2)$), and that there exists some $s > \alpha -1$ such that $\omega^s$ satisfies the condition $(\alpha_0)$. Let us observe that we can take $s'<s$ such that $1>s'> \alpha -1$, and it is then easy to show that $\omega^{s'}$ satisfies the conditions $(\omega_5)$ and $(\alpha_0)$.
        %\red{, so that Theorem~\ref{theoremWeightFunctionWideSector} can be applied to obtain the equivalence of (a) through (f).}

     \item[(ii)] If $\alpha\geq 2$, suppose there exists $s$ according to the assumptions in the theorem. Then, we will have $s>1$, and since $\omega^s$ satisfies the condition $(\omega_5)$, we can check immediately that $\omega$ has $(\omega_2)$.
         %\red{ So, again all the conditions in the theorem are equivalent.}
\end{itemize}

\end{remark}

\section{Examples}\label{sect.Examples}
{In this section, we apply the previous results to some well-known examples of ultraholomorphic classes. Let us fix $\alpha>0$.

\subsection{Gevrey-related classes}

Consider the sequence $\overline{\G}^\beta:=(j^{j\beta})_{j\in\NN_0}$ of index $\beta\in\RR$. Note that this sequence has the \hyperlink{rai}{$(\on{rai})$} property if and only if $\beta\geq1$. We are going to study the stability of the class $\mathcal{A}_{\{\overline{\G}^\beta\}}(S_{\alpha})$ in terms of the values of $\alpha$ and $\beta$. Let us distinguish some cases:

\begin{enumerate}
	\item[(a)] Let $\alpha\in(0,1]$:
				\begin{enumerate}
					\item[(i)] If $\beta<\alpha-1$ then  $\lim_{j\rightarrow+\infty}(j^{(1-\alpha)j}j^{j\beta})^{1/j}=0$, and therefore the class is stable because it is trivial, i.e., it only contains constant functions (see Remark \ref{liminf-remark}).
					\item[(ii)] If $\beta\in(\alpha-1,1)$ Corollary \ref{coroCarlemanClassNarrowSector}, together with the fact that $\overline{\G}^\beta$ has not the \hyperlink{rai}{$(\on{rai})$} property, ensure that the class is non stable.
					\item[(iii)] If $\beta=\alpha-1$, the sequence $\M^{\alpha}$ is $\overline{\G}^\beta$, which does not satisfy (rai). So, by Remark \ref{ramimatrixdc} and Corollary \ref{coroCarlemanClassNarrowSector} the class is not stable.
					\item[(iv)] If $\beta\geq1$ we deduce from the Corollary \ref{coroCarlemanClassNarrowSector} that the class is stable.
				\end{enumerate}
	\item[(b)] Let $\alpha>1$:
				\begin{enumerate}
					\item[(i)] If $\beta\leq\alpha-1$ then the $\liminf_{j\rightarrow+\infty}(j^{(1-\alpha)j}j^{j\beta})^{1/j}<\infty$, and therefore the class is stable because it only contains constant functions (see Remark \ref{liminf-remark}).
					\item[(ii)] If $\beta>\alpha-1$, we have stability provided that $\beta\geq1$, thanks to the Corollary \ref{coroCarlemanClassWideSector}.
%					\item[(iii)] If $\alpha\geq 2$ and $\beta>\alpha-1$, then $\overline{\G}^\beta$ has the \hyperlink{rai}{$(\on{rai})$} property, and the Corollary \ref{coroCarlemanClassWideSector} ensures that the class is stable.
				\end{enumerate}
\end{enumerate}
We include a graphic in order to see the stability (resp. non stability) regions:
}
\begin{figure}[!h]
	\centering
	\begin{tikzpicture}[line cap=round,line join=round,>=triangle 45,x=1.0cm,y=1.0cm]
		\begin{axis}[
			x=1.0cm,y=1.0cm,
			axis lines=middle,
			xtick={-0,1,...,3},
			ytick={-2,-1,...,2},
			xlabel={$\alpha$},
			ylabel={$\beta$},
			xlabel style={below right},
			ylabel style={above left},
			xmin=-0.5,
			xmax=3.5,
			ymin=-2.0,
			ymax=3.2]
			\draw[color=ududff] (2.5,1.9) node[anchor=center,rotate=45] {$\beta=\alpha-1$};
			\clip(-0.1,-2.) rectangle (3.3,2.9);
			\draw [line width=2.8pt,color=ududff] (0.,1.)-- (2.,1.);
			\draw [line width=2.8pt,color=ududff] (2.,1.)-- (1.,0.);
			%\draw [line width=2.8pt,dash pattern=on 1pt off 1pt,color=ffqqqq] (1.,0.)-- (0.,-1.);
			\draw [line width=2.8pt,color=ududff] (0.,-2.)-- (0.,-1.);
			\draw [line width=2.8pt,color=ududff] (0.,1.) -- (0.,3.0);
			%\draw [->,line width=1.0pt,color=ududff] (0.,1.) -- (0.,2.1673608209628754);
			\draw [line width=2.pt,color=ududff,domain=1.4428914522128231:3.4] plot(\x,{(-0.5571085477871769--0.5571085477871769*\x)/0.5571085477871769});
			\draw[line width=0.pt,dash pattern=on 1pt off 1pt,color=qqqqff,fill=qqqqff,fill opacity=0.2](3.4,2.45)--(0,-1)--(0,-2.1)--(3.4,-2.1)--(3.4,2.45); %Trapecio 1
			 %%%%%%%%%%%%%%%%%%%%%%%%%%%%%%%%%%%%%%%%%%%%%%%%%%%%%%%%%%%%%%%%%%%%%%%%%%%%%%%%%%%%%%%%%%%%%%%%%%%%%%%%%%%%%%%%%%%%%%%%%%%%%%%%%%%%%%%%
			%\draw[line width=0.pt,dash pattern=on 1pt off 1pt,color=qqqqff,fill=qqqqff,fill opacity=0.20000000298023224](3.3446367832674433,2.3)--(2.6495045587295146E-8,-0.9999999383648361)--(2.6495045587295146E-8,-1.6686155136190581)--(3.6111876344354417,-1.6686155136190581)--(3.6111876344354417,2.17);
			 %%%%%%%%%%%%%%%%%%%%%%%%%%%%%%%%%%%%%%%%%%%%%%%%%%%%%%%%%%%%%%%%%%%%%%%%%%%%%%%%%%%%%%%%%%%%%%%%%%%%%%%%%%%%%%%%%%%%%%%%%%%%%%%%%%%%%%%%
			\draw[line width=0.pt,dash pattern=on 1pt off 1pt,color=qqqqff,fill=qqqqff,fill opacity=0.1](0,3)--(0,1.000000069248137)--(2.000000034108022,1.000000069248137)--(3.4,2.45)--(3.4,3);%Trapecio 2
			 %%%%%%%%%%%%%%%%%%%%%%%%%%%%%%%%%%%%%%%%%%%%%%%%%%%%%%%%%%%%%%%%%%%%%%%%%%%%%%%%%%%%%%%%%%%%%%%%%%%%%%%%%%%%%%%%%%%%%%%%%%%%%%%%%%%%
			%\draw[line width=0.pt,dash pattern=on 1pt off 1pt,color=qqqqff,fill=qqqqff,fill opacity=0.20000000298023224](2.6495045587295146E-8,2.17)--(2.6495045587295146E-8,1.000000069248137)--(2.000000034108022,1.000000069248137)--(3.3446367832674433,2.17);
			 %%%%%%%%%%%%%%%%%%%%%%%%%%%%%%%%%%%%%%%%%%%%%%%%%%%%%%%%%%%%%%%%%%%%%%%%%%%%%%%%%%%%%%%%%%%%%%%%%%%%%%%%%%%%%%%%%%%%%%%%%%%%%%%%%%%%
			\draw[line width=0.pt,dash pattern=on 1pt off 1pt,color=ffqqqq,fill=ffqqqq,fill opacity=0.25](2.6495045587295146E-8,1.000000069248137)--(2.6495045587295146E-8,-0.9999999383648361)--(2.000000034108022,1.000000069248137);%Triángulo
			 %%%%%%%%%%%%%%%%%%%%%%%%%%%%%%%%%%%%%%%%%%%%%%%%%%%%%%%%%%%%%%%%%%%%%%%%%%%%%%%%%%%%%%%%%%%%%%%%%%%%%%%%%%%%%%%%%%%%%%%%%%%%%%%%%%%%%
			%\draw[line width=0.pt,dash pattern=on 1pt off 1pt,color=ffqqqq,fill=ffqqqq,fill opacity=0.25](2.6495045587295146E-8,1.000000069248137)--(2.6495045587295146E-8,-0.9999999383648361)--(2.000000034108022,1.000000069248137);
			 %%%%%%%%%%%%%%%%%%%%%%%%%%%%%%%%%%%%%%%%%%%%%%%%%%%%%%%%%%%%%%%%%%%%%%%%%%%%%%%%%%%%%%%%%%%%%%%%%%%%%%%%%%%%%%%%%%%%%%%%%%%%%%%%%%%%%
			\draw[line width=0.pt,dash pattern=on 1pt off 0pt,color=green,fill=green,fill opacity=0.07](3.4,2.4)--(0,-1)--(0,-2.1)--(3.4,-2.1)--(3.4,2.4); %Trapecio 1 verde
			 %%%%%%%%%%%%%%%%%%%%%%%%%%%%%%%%%%%%%%%%%%%%%%%%%%%%%%%%%%%%%%%%%%%%%%%%%%%%%%%%%%%%%%%%%%%%%%%%%%%%%%%%%%%%%%%%%%%%%%%%%%%%%%%%%%%%%%%%
			%
			%\draw [color=ffqqqq](0.8,0.7) node[anchor=center,rotate=20] {\tiny Non stable};
			\draw [color=ffqqqq](0.9,0.25) node[anchor=center,rotate=45] {\small \textbf{Non stable}};
			%\draw [color=ududff](1.5,-1) node[anchor=center] {\small \textbf{Stable}};
			\draw [color=ududff](1.5,2.3) node[anchor=center] {\small \textbf{Stable}};
			\draw [color=verdeosc](2,.4) node[anchor=center,rotate=45] {\small \textbf{Trivial class}};
			%%%%%%%%%%%%%%% Segmentos %%%%%%%%%%%%%%%%%%%%%%%%%%%%%%%%%%%%%
			\draw [line width=2.8pt,color=ududff] (0.,1.)-- (2.,1.);
			\draw [line width=2.8pt,color=ududff] (2.,1.)-- (1.,0.);
			\draw [line width=2.8pt,dash pattern=on 1pt off 5pt,color=ffqqqq] (1.,0.)-- (0.,-1.);
			\draw [line width=2.8pt,color=ududff] (0.,-2.)-- (0.,-1.);
			\draw [line width=2.8pt,color=ududff] (0.,1.) -- (0.,3.0);
			\draw [line width=2.8pt,dash pattern=on 1pt off 5pt,color=ffqqqq] (0.,-1.)-- (0.,1.);% Línea discontinua
			%\draw [->,line width=1.0pt,color=ududff] (0.,1.) -- (0.,2.1673608209628754); %Vector o Segmento flechado
			\draw [line width=2.pt,color=ududff,domain=1.4428914522128231:3.4] plot(\x,{(-0.5571085477871769--0.5571085477871769*\x)/0.5571085477871769}); %Semirrecta \beta=\alpha-1
			%%%%%%%%%%%%%%%%%%%%%%%%%%%% Puntos %%%%%%%%%%%%%%%%%%%%%%%%%%%%%%%%%%%
			\begin{scriptsize}
				\draw [fill=xdxdff] (0.,1.) circle (2.5pt);
				\draw [fill=ffqqqq] (0.,-1.) circle (2.5pt);
				\draw [fill=ffqqqq] (1.,0.) circle (2.5pt);
				%\draw [color=xdxdff] (1.,0.) circle (2.5pt);% Punto hueco
				\draw [fill=ududff] (2.,1.) circle (2.5pt);
				%\draw [fill=xdxdff] (0.,-2.) circle (2.5pt);
				%\draw [fill=xdxdff] (0.,2.1673608209628754) circle (2.5pt);
				%\draw[color=ududff] (1.5,1.6) node {$\beta=\alpha-1$};
				%%\draw[color=ududff] (2.5,1.9) node[anchor=center,rotate=45] {$\beta=\alpha-1$};
			\end{scriptsize}
		\end{axis}
	\end{tikzpicture}
\end{figure}
{
We consider now a second example. Let us fix $\alpha>1$, take some $\beta>\alpha$ and consider the weight matrix $\mathcal{L}^{(\beta)}=\{\overline{\G}^{\beta-\frac{1}{p+1}}: p>0\}$. Note that the ultraholomorphic class associated with $\mathcal{L}^{(\beta)}$ is strictly smaller than the class associated with the constant matrix $\mathcal{G}^{\beta}=\{\overline{\G}^{\beta}: p>0\}$. Under these assumptions, let us observe that $\overline{\G}^{\beta-\frac{1}{p+1}}$ is an (lc) sequence for all $p>0$. Then Theorem \ref{theorem1siddiqi.} guarantees that the class $\mathcal{A}_{\{\mathcal{L}^{(\beta)}\}}(S_{\alpha})$ is stable, thanks to the fact that $\beta-\frac{1}{p+1}>1$ for large $p$, and we can ensure that the corresponding matrix has \hyperlink{R-rai}{$(\mathcal{M}_{\{\on{rai}\}})$}.

\subsection{q-Gevrey case}

In this subsection, we will work, for $q>1$, with the {\itshape q-Gevrey sequence}, i.e $\M_q=(q^{j^2})_{j\geq 0}$.
First, thanks to the fact that the sequence $\M_q$ has (lc) and (dc), and moreover $\check{M}_q$ is also (lc), we can easily prove the stability properties for the class $\mathcal{A}_{\{\M_q\}}(S_{\alpha})$. For $\alpha\in (0,1]$, the Corollary \ref{coroCarlemanClassNarrowSector} ensures that the class $\mathcal{A}_{\{\M_q\}}(S_{\alpha})$ is stable. On the other hand, for $\alpha>1$ and for any $\beta>\alpha$ the sequence $\overline{\G}^{1-\beta}\M_q$ is equivalent to an (lc) sequence, because the gamma index of $\M_q$ is  infinity. So, the Corollary \ref{coroCarlemanClassWideSector} again ensures the stability.

Now, we want to study the stability properties for the class $\mathcal{A}_{\{\omega_{\M_q}\}}(S_{\alpha})$. For this purpose, let us observe that we can estimate the normalized weight function $ \omega_{\M_q}$,
$$\omega_{\M_q}(t)= \sup_{j\in\mathbb{N}_0} \ln\left(\frac{t^{j}}{q^{j^2}}\right)=\sup_{j\in\mathbb{N}_0} (j\ln(t)-j^2\ln(q)),\quad t>1.
$$
Obviously, $\omega_{\M_q}(t)$ is bounded above by the supremum of $x\ln(t)-x^2\ln(q)$ when $x$ runs over $(0,\infty)$, which is easily obtained by elementary calculus and occurs at the point
$$
\left(\frac{\ln(t)}{2\ln(q)},\frac{\ln^2(t)}{4\ln(q)}\right).
$$
In particular, it is easy to check that $\omega(t):=\ln^2(t)/(4\ln(q))$ verifies (after normalization in the interval $[0,1]$) that $\omega\in\hyperlink{omset1}{\mathcal{W}}$, $\omega$ has $(\omega_5)$ (and therefore $(\omega_2)$) and $\omega\hyperlink{sim}{\sim}\omega_{\M_q}$, so the corresponding matrices $\mathcal{M}_\omega$ and $\mathcal{M}_{\omega_{\M_q}}$ are $R$-equivalent. In order to compute the matrix associated with $\omega$, the Legendre-Fenchel-Young-conjugate of $\varphi_{\omega}$ is
$$
\varphi^{*}_{\omega}(x):=\sup_{y\ge 0}\{x y-\omega(\exp(y))\}=x^2\ln(q)=\ln(q^{x^2}),\qquad x\ge 0.
$$
So, we have that
$$
W^{(\ell)}_j=\exp(\frac{1}{\ell}\varphi^{*}_{\omega}(\ell_j))=q^{\ell j^2},\qquad j\geq0, \text{ and therefore}\qquad \W^{(\ell)}=(q^{\ell j^2})_{j\geq 0},\qquad \ell>0.
$$
Note that each sequence $\W^{(\ell)}$ is (lc), (dc) and has the property (rai) for all $\ell>0$, in this situation Theorem \ref{theorem1siddiqcase1} ensures that the class $\mathcal{A}_{\{\omega\}}(S_{\alpha})$ (resp. $\mathcal{A}_{\{\omega_{\M_q}\}}(S_{\alpha})$) is stable for $\alpha\in(0,1]$. On the other hand, note that $\gamma(\omega)=\infty$, since $\gamma(\omega)\geq \gamma(\W^{(\ell)})$ for all $\ell>0$ (see Subsection~\ref{ssect.IndexGammaOmega}) and $\gamma(\W^{(\ell)})$ is also infinity. In this case, Remark \ref{gammaindex-omega} ensures that we can apply Theorem \ref{theoremWeightFunctionWideSector} in order to deduce that the class $\mathcal{A}_{\{\omega\}}(S_{\alpha})$ (resp. $\mathcal{A}_{\{\omega_{\M_q}\}}(S_{\alpha})$) is stable for $\alpha>1$.

}

\vskip.5cm
\noindent\textbf{Acknowledgements}: The first three authors are partially supported by the Spanish Ministry of Science and Innovation under the project PID2019-105621GB-I00. The fourth author is supported by FWF-Project P33417-N.

\bibliographystyle{plain}

%\vskip.4cm
%\noindent\textbf{Statements and declarations}:\\
%\vskip.1cm\noindent \textbf{Funding} \ The first three authors are partially supported by the Spanish Ministry of Science and Innovation under the project PID2019-105621GB-I00. The fourth author is supported by FWF-Project P33417-N.\par
%
%\vskip.2cm\noindent \textbf{Competing interests} \ The authors have no relevant financial or non-financial interests to disclose.\par
%
%\vskip.2cm\noindent \textbf{Author contributions} \ All authors contributed to the study and the preparation of previous versions of the manuscript. All authors read and approved the final manuscript.\par
%
%\vskip.2cm\noindent \textbf{Data availability} \ Data sharing not applicable to this article as no datasets were generated or analysed during the current study.

\vskip.5cm
\noindent\textbf{Affiliations}:\\
\noindent Javier~Jim\'{e}nez-Garrido:\\
Departamento de Matem\'aticas, Estad{\'\i}stica y Computaci\'on\\
Universidad de Cantabria\\
Avda. de los Castros, s/n, 39005 Santander, Spain\\
Instituto de Investigaci\'on en Matem\'aticas IMUVA, Universidad de Va\-lla\-do\-lid\\
ORCID: 0000-0003-3579-486X\\
E-mail: jesusjavier.jimenez@unican.es\\

\vskip.1cm
\noindent
Ignacio Miguel-Cantero:\\
Departamento de \'Algebra, An\'alisis Matem\'atico, Geometr{\'\i}a y Topolog{\'\i}a\\
Universidad de Va\-lla\-do\-lid\\
Facultad de Ciencias, Paseo de Bel\'en 7, 47011 Valladolid, Spain.\\
Instituto de Investigaci\'on en Matem\'aticas IMUVA\\
ORCID: 0000-0001-5270-0971\\
E-mail: ignacio.miguel@uva.es\\

\vskip.1cm
\noindent Javier~Sanz:\\
Departamento de \'Algebra, An\'alisis Matem\'atico, Geometr{\'\i}a y Topolog{\'\i}a\\
Universidad de Va\-lla\-do\-lid\\
Facultad de Ciencias, Paseo de Bel\'en 7, 47011 Valladolid, Spain.\\
Instituto de Investigaci\'on en Matem\'aticas IMUVA\\
ORCID: 0000-0001-7338-4971\\
E-mail: javier.sanz.gil@uva.es\\

\vskip.1cm
\noindent Gerhard~Schindl:\\
Fakult\"at f\"ur Mathematik, Universit\"at Wien,
Oskar-Morgenstern-Platz~1, A-1090 Wien, Austria.\\
ORCID: 0000-0003-2192-9110\\
E-mail: gerhard.schindl@univie.ac.at
\end{document}